\documentclass[a4paper,12pt,leqno]{article}

\usepackage{geometry}
\usepackage{tikz-cd}\usepackage{bussproofs}
\usepackage{amssymb,amsmath}
\usepackage{proof}
\usepackage{psvectorian}
\usepackage{CJKutf8}
\usepackage[utf8]{inputenc}
\usepackage{graphicx}
\usepackage{graphicx}
\graphicspath{ {images/} }
\usepackage{qtree}
\usepackage{mathtools}
\usepackage{hyperref}
\usepackage{listings}
\usepackage{indentfirst}
\usepackage{titlesec}
\numberwithin{equation}{section}

\newtheorem{defn}[equation]{Definition}

\newtheorem{rem}[equation]{Remark}
\newtheorem{exm}[equation]{Example}

\newtheorem{theorem}[equation]{Theorem}
\newtheorem{notat}[equation]{Notation}
\newtheorem{newpar}[equation]{}

\newtheorem{xdefn}{Definition.}
\newtheorem{xproposition}{Proposition.}
\newtheorem{xcorollary}{Corollary.}
\newtheorem{xrem}{Remark.}
\newtheorem{xexm}{Example.}
\newtheorem{xlemma}{Lemma.}
\newtheorem{xtheorem}{Theorem.}
\newtheorem{xnotat}{Notation.}
\newtheorem{xnewpar}{\it}
\newtheorem{xproof}{{\it Proof. }}
\newtheorem{xproofof}{{\it Proof}}

\newenvironment{proof}{\begin{xproof}\em}{\end{xproof}}

\newenvironment{newparagraph*}[1]{\begin{xnewpar}\hspace*{-1.5mm}{#1}. \rm}{\end{xnewpar}}

\newenvironment{definition*}{\begin{xdefn}\em}{\end{xdefn}}
\newenvironment{remark*}{\begin{xrem}\em}{\end{xrem}}
\newenvironment{example*}{\begin{xexm}\em}{\end{xexm}}
\newenvironment{notation*}{\begin{xnotat}\em}{\end{xnotat}}
\newenvironment{proposition*}{\begin{xproposition}}{\end{xproposition}}
\newenvironment{corollary*}{\begin{xcorollary}}{\end{xcorollary}}
\newenvironment{lemma*}{\begin{xlemma}}{\end{xlemma}}
\newenvironment{theorem*}{\begin{xtheorem}}{\end{xtheorem}}

\titleformat*{\section}{\large\bfseries}

\begin{document}

\title{Aristotle's Second-Order Logic and Natural Deduction}
\author{Clarence Protin\footnote{Centro de Filosofia da Universidade de Lisboa.}}

\date{18th of January of 2026}
\maketitle
\begin{abstract}
This paper has two goals. The first goal is to show how an extension of second-order logic is a natural framework to  formalize portions of Aristotle's \emph{Topics} and to bring to the foreground the logical, linguistic and philosophical interest of this  work, showing in particular that we are in the presence of a richly intensional and modal conception of logic.  Aristotelian logic  and its related traditions in antiquity are often held to have been equivalent to monadic predicate logic and as such inadequate to formalize mathematics as well as scientific and philosophical discourse in general. The second goal of this paper is to argue that on the contrary the logical theories of Aristotle (which we argue correspond to a variant of natural deduction)  and ancient authors  such as Galen and Boethius were in fact  quite sufficient to account for the logically complex expressions and  reasoning involving multiple generality fundamental to the aforementioned disciplines.

\end{abstract}


\section{Introduction}

This paper will be organized as follows.  Sections 2-6 are dedicated to the formalization and discussion of Aristotle'\emph{Topics} using an extension of second-order logic. In Sections 7-9 we argue that the logical theory of Aristotle and some ancient authors was fundamentally
equivalent to a variant of natural deduction which was applied to formulas of any complexity.
In Section 2 we present the language and rules of an extension of classical second-order logical (ESOL) which will serve as the foundation for Sections 3-6.  In Section 3 we  discuss briefly some aspects of metalogic of the \emph{Analytics}  and the formal structure of ancient logic which will be followed by a  fuller discussion in Section 9.   In Section 4 we study some key aspects of Aristotle's \emph{Topics} and argue that ESOL is a good choice for representing the underlying logic. Here we require a richer collection of second-order predicates (expressing the predicabilia) as compared to the \emph{Analytics}; we do not have a  'extensionalist' reduction of the main second-order predicates to formulas involving only first-order predication.
 A \emph{topic} is a sentence of ESOL in general universally quantified over most of the monadic predicate variables occurring in the sentence. These universally quantified monadic predicate variables are then supposed to be instantiated according to the dialectical needs at hand\footnote{The definition of topic (locus) given by Boethius  in De topicis differentiis agrees with this. The definition 'Locus autem sedes est argumenti, uel id unde ad propositam quaestionem conueniens trahitur argumentum' suggests directly a second-order universal quantified sentence  which can be instantiated to yield an inference scheme to be used according to the 'predicates'  at hand, the proposita questio.}. 
 Although the structure of the \emph{Topics} consists mainly in a rather tedious enumeration of topics, the most purely logical among them can be organized into groups wherein all the topics in a group are logical consequences of a certain paradigmatic topic. Often Aristotle indicates explicitly that for two closely related topics (mentioned in succession) one is the consequence of the other. In Section 5 we look at book V which deals with property (\emph{idion}) and discuss some of the challenges facing a formalization of the topics included therein as well as analyzing Aristotle's treatment of superlatives in book VII, 152a.
  In Section 6 we end with a discussion about how the system of debate for which the \emph{Topics} was written might be reconstructed and interpreted according to ESOL.    In order to provide more context for Sections 4-5, in appendix 1 we list a formalization of an important subset of the topics of book IV, focused on the relation between genus and species. 
In Section 7 we set up a variant of ESOL in which all quantifiers are bounded by some formula and the elimination rules for $\exists$ are replaced by rules involving an indefinite version of the Peano operator. The resulting system is arguably closer to the logical and inferential mechanisms of natural language and  provides
a good foundation to argue for the main thesis of the second part of the paper. A common view is that whereas  the deductive mechanisms and logical expressivity of ancient Greek mathematics corresponds to  the full power of modern first or higher order predicate logic, it was nonetheless the case that Aristotle's theory of deduction (identified with the theory of the  syllogism of the Prior Analytics ) captures only monadic predicate logic and was thus radically inadequate for the logical structure and deductive mechanism of mathematics as well as that of scientific and philosophical discourse in general. 
  In Sections 8 and 9 we challenge the above view and argue that the logical theories of Aristotle and some authors of late antiquity (Galen and Boethius) were  in fact sufficient to account for the logically complex expressions and deductive processes of ancient mathematics, science and philosophy.
  To this end in Section 8 we show first that Aristotle not only makes use of the four quantifier rules for ESOL$_D$ adequate for reasoning about quantified sentences but does so in a philosophical reflected way, that is, the quantifier rules are part of his theory of logic. We also discuss his use of reductio ad absurdum and the rules related to implication and inference.  We then show that these rules are deployed by Aristotle for sentences containing multiple quantifiers over relations and furthermore present  evidence  for such a deployment in Galen.  Section 9 is dedicated to discussing evidence for a formal theory of multiple generality in Boethius. For Aristotle we refer to the edition of the Greek text in \cite{topiques1} and \cite{topiques2} for the Topics, \cite{Ros} for the Prior and Posterior Analytics and \cite{Ros2} for the Physics.

 \section{Extended Second-Order Logic}

 The language of extended second-order logic (ESOL) is as follows. We are given a collection of individual variables $x,y,z,...$ and individual constants $a,b,c,...$ together with a collection of $n$-ary predicate term variables $X^{(n)}, Y^{(n)}, Z^{(n)},...$ (we omit the superscript if the context is clear) for $n\geq 0$. Furthermore we are given a collection of predicate term symbols $\alpha^{(\tau)}, \beta^{(\tau)}, \gamma^{(\tau)},...$ where $\tau$ consists of a (possibly empty) list $(n_1,...,n_m)$ of integers $n_i \geq -1$ called a signature,  where we likewise omit the superscript if the context is clear. It is also convenient to consider that we have a collection of function symbols $f^{(\tau,n)}, g^{(\tau,n)},...$ where $\tau$ is a signature and $n\geq -1$. For $n = -1$ these are called first-order functional symbols. Otherwise they are called second-order function symbols.
 
 We define formula and $n$-predicate term and individual term by mutual recursion.
 
 For convenience we use small-case Roman letters to represent  individual terms, variables and constants and upper-case Roman letters to represent  $n$-predicate terms and variables. We use gothic letters to represent terms in general.  We use upper-case Greek letters to represent formulas.
 
 If $\Phi_1, \Phi_2$ are formulas then so are $(\Phi_1 \vee \Phi_2), (\Phi_1 \wedge \Phi_2), \neg \Phi_1, \forall x \Phi_1, \exists x \Phi_1, \forall X \Phi_1, \exists X \Phi_1$ where $x$ is an individual variable and $X$ is a predicate term variable.

 If $\alpha^{(\tau)}$ is a predicate term symbol with $\tau = (n_1,...,n_m)$ non-empty then $\alpha \mathfrak{t}_1....\mathfrak{t}_m$ is a formula where if $n_i \geq 0$ then $\mathfrak{t}_i$ is a $n_i$-predicate term and if $n_i = -1$ then $\mathfrak{t}_i$ is an individual term. If $\tau$ is empty then $\alpha$ is considered a $0$-predicate term and we write the corresponding formula as $\alpha()$.

 If $\Phi$ is a formula and $x_1,...,x_n$ is a list of individual variables (possible empty) then $\lambda x_1...x_n \Phi$ is a $n$-predicate term (for the empty list we put $n= 0$ and write $\lambda \Phi$).
 
 If $X^{(n)}$ is a predicate term variable then it is a $n$-predicate term.
 If $\alpha^{(\tau)}$ is a predicate term symbol such that $\tau= (k_1,...,k_n)$ with $k_i = -1$ then $\alpha^{(\tau)}$ is a $n$-predicate term\footnote{we can also call $\alpha$ in this case a $n$-predicate \emph{constant}.}
 
 If $T$ is a $n$-predicate term  and $t_1,...t_n$ are individual terms then $Tt_1...t_n$ is a formula. For the case $n = 0$ we write the formula $T()$.

 An individual constant or variable is an individual term.

 If $f^{(\tau,n)}$ is a function symbol then $f\mathfrak{t}_1...\mathfrak{t}_m$  is an individual term (for $n= -1$) or  $n$-predicate term (for $n\geq 0$) where $\mathfrak{t}_1...\mathfrak{t}_m$ satisfy the condition relative to $\tau$ as above.

 The  standard natural deduction rules consist of the standard rules for the first-order classical predicate calculus:
 
 The rules are as follows:

 \begin{center}
 	\AxiomC{ $\Phi$  \quad $\Psi$ }
 	\RightLabel{$\& I$}
 	\UnaryInfC{$\Phi\,\&\,\Psi$}
 	
 	\DisplayProof
 	\AxiomC{ $\Phi\,\&\, \Psi$ }
 	\RightLabel{$\& E_R$}
 	\UnaryInfC{$\Phi$}
 	\DisplayProof
 	\AxiomC{ $\Phi\,\&\, \Psi$ }
 	\RightLabel{$\& E_L$}
 	\UnaryInfC{$\Psi$}
 	\DisplayProof
 \end{center}
 
 \begin{center}
 	\AxiomC{ $\Phi$ }
 	\RightLabel{$\vee I_R$}
 	\UnaryInfC{$\Phi\,\vee\,\Psi$}
 	\DisplayProof 
 	\AxiomC{ $\Psi$ }
 	\RightLabel{$\vee I_L$}
 	\UnaryInfC{$\Phi\,\vee\,\Psi$}
 	\DisplayProof 
 	\AxiomC{$\Phi \lor \Psi$}
 	\AxiomC{[$\Phi$]}
 	\noLine
 	\UnaryInfC{$\Upsilon$}
 	\AxiomC{[$\Psi$]}
 	\noLine
 	\UnaryInfC{$\Upsilon$}
 	\RightLabel{$\vee E$}
 	\TrinaryInfC{$\Upsilon$}
 	\DisplayProof
 \end{center}
 \begin{center}
 	\AxiomC{($\Phi$)}
 	\noLine
 	\UnaryInfC{$\Psi$}
 	\RightLabel{$\rightarrow$I}
 	\UnaryInfC{$\Phi\rightarrow \Psi$}
 	\DisplayProof
 	\AxiomC{$\Phi$}
 	\AxiomC{$\Phi\rightarrow \Psi$}
 	\RightLabel{$\rightarrow$E}
 	\BinaryInfC{$\Psi$}
 	\DisplayProof
 \end{center}
 \begin{center}
 	\AxiomC{$\bot$}
 	\RightLabel{$\bot_i$}
 	\UnaryInfC{$\Phi$}
 	\DisplayProof
 	\AxiomC{$(\sim \Phi)$}
 	\noLine
 	\UnaryInfC{$\bot$}
 	\RightLabel{$\bot_c$}
 	\UnaryInfC{$\Phi$}
 	\DisplayProof
 \end{center}
 
 the four first-order quantifier rules:
 
 \begin{center}
 	\AxiomC{$\Phi$}
 	\RightLabel{$\forall^1 I$}
 	\UnaryInfC{$\forall x\Phi^y_x$}
 	\DisplayProof
 	\AxiomC{$\forall x\Phi$}
 	\RightLabel{$\forall^1 E$}
 	\UnaryInfC{$\Phi^x_t$}
 	\DisplayProof
 \end{center}
 
 \begin{center}
 	\AxiomC{$\Phi^x_t$}
 	\RightLabel{$\exists^1 I$}
 	\UnaryInfC{$\exists x\Phi$}
 	\DisplayProof
 	\AxiomC{$\exists x\Phi$}
 	\AxiomC{$(\Phi^x_y)$}
 	\noLine
 	\UnaryInfC{$\Psi$}
 	\RightLabel{$\exists^1 E$}
 	\BinaryInfC{$\Psi$}
 	\DisplayProof
 \end{center}

 where in $\forall^1 I$ $y$ must not occur free in any assumption on which $\Phi$ depends and in $\exists^1 I$ $y$ must not occur free in $\Psi$ or in $\exists x \Phi$ or on any hypothesis on which $\Psi$ depends except for $\Phi^x_y$.
 Here $\Phi^t_s$ denotes the result of substituting $t$ by $s$ in $\Phi$\footnote{Here we need not follow Prawitz's terminological distinction between variables and parameters which correspond to bound and free variables.}.

 We have the analogous second-order quantifier rules

 \begin{center}
 	\AxiomC{$\Phi$}
 	\RightLabel{$\forall^2 I$}
 	\UnaryInfC{$\forall X^n\Phi^{Y^n}_{X^n}$}
 	\DisplayProof
 	\AxiomC{$\forall X^n\Phi$}
 	\RightLabel{$\forall^2 E$}
 	\UnaryInfC{$\Phi^{X^n}_{T^n}$}
 	\DisplayProof
 \end{center}
 
 \begin{center}
 	\AxiomC{$\Phi^{X^n}_{T^n}$}
 	\RightLabel{$\exists^2 I$}
 	\UnaryInfC{$\exists X^n\Phi$}
 	\DisplayProof
 	\AxiomC{$\exists X^n\Phi$}
 	\AxiomC{$(\Phi^{X^n}_{T^n})$}
 	\noLine
 	\UnaryInfC{$\Psi$}
 	\RightLabel{$\exists^2 E$}
 	\BinaryInfC{$\Psi$}
 	\DisplayProof
 \end{center}
 
 with analogous provisos.
 
 Finally there are introduction and elimination rule for $\lambda$-predicate terms for which the reader may consult \cite{pra}. In Prawitz\cite{pra} versions 1 and 2 are proof-theoretically equivalent (for formulas without $\lambda$), the difference being that in version 1 the rules $\lambda I$ and $\lambda E$ are incorporated into second-order quantifier rules\footnote{Why do we consider this system as being an extension of second-order logic rather than a fragment of third-order logic ? We are of the view that what defines the 'order' of a logic is what is quantified over and what kind of corresponding variable there is. First-order logic quantifies over individuals and has individual variables. Second-order logic quantifies over entities corresponding to $n$-ary predicates with individual arguments and has the corresponding $n$-ary predicate term variables. Third-order logic thus would have to quantify over entities corresponding to predicates which can involve both individual arguments and $n$-ary predicates of individuals arguments and furthermore it would require variables to represent these last. Such a quantification and such variables are lacking in this system although the predicate term symbols $\alpha^{(\tau)}$ can be considered in some sense as being 'third-order' when involving $n$-predicate term arguments.}.  In this paper we only consider classical logic (i.e. we always have the rule $\bot_c$).

 \section{Metalogic of the Analytics}

 \subsection{Logic and language}

 There have been many attempts to formalize Aristotle's syllogistic expounded in the \emph{Prior Analytics} in both its strictly assertoric and modal variants (see \cite{pro} for references). Some approaches involve monadic (modal or otherwise) first-order logic and other approaches use variants (modal of otherwise) of propositional logic. Additionally many authors assume that Aristotle reasons much  like a modern logician by constructing models to refute the validity of certain syllogistic figures. 
 The above approaches usually set up a minimal set of axioms together with a deductive system and attempt to show the agreement (at in least for most cases)  between the derivable syllogisms and those considered valid by Aristotle (how Aristotle himself derived or justifies the validity or non-validity playing a secondary role).  Our concern in this Section is different. We are interested in the structure of the proofs in the \emph{Analytics} themselves, how Aristotle argues for the validity or non-validity of a certain figure or propositional conversion - an approach reminiscent of the modern division between metalogic and object logic.  
 A noteworthy characteristic of the \emph{Analytics} is the presence of variable letters ranging over Aristotelian \emph{terms} which (so we argue) usually correspond to monadic predicate terms in second-order logic. There is a distinguished class of entities for which it is meaningful to ask if such entities fall or do not fall under such monadic predicate terms - these we call \emph{individuals}. Thus we can ask if a monadic predicate term is predicated or not of a given individual. 
 
 It is patent that a formalization of Aristotle's reasoning must be a term-logic, that is, terms must be the basic building blocks - the most basic division being between individual terms (corresponding to 'first substances') and non-individual terms. Sentences (note that the \emph{heis logos}, \emph{De Int.} 17a16-17 cannot be identified with 'atomic' sentences in the modern sense as they may contain connectives)  are then built up from combinations of terms in a specified manner. Predication - and in Aristotle there are many distinct notions of predication hiding under the verb $\emph{esti}$ - corresponds to combinations of terms tagged according to the kind of predication being considered. Thus if \textcircled{1} corresponded to predicating a non-individual term $A$ of an individual term $a$ (let us say according to 'essence') and \textcircled{2} corresponded to predicating a non-individual term $B$ of another $A$ according to 'essence' then grammatically (for SOV languages like Ancient Greek) and notationally the most natural representation would be aA\textcircled{1} and AB\textcircled{2}. Such a representation is necessarily if we wish to formalize correctly  Aristotle's essentially intensional logic: the predication of definition and the predication of property (\emph{idion}) both imply that the terms have equal 'extension' (in the sense of individuals falling under them) but the types of predication are different as are the terms: 'being able to laugh' is not the same term as 'rational animal' but both are co-extensional with 'man'.  Thus we have grounds for holding that an extension of second-order logic  (in which \textcircled{2} corresponds to a predicate term symbol of signature $(1,1)$ and \textcircled{1} to simple concatenation in reverse order $Aa$) is an adequate minimal formalization for the logical-argumentative apparatus in the \emph{Organon}.
 What we wrote above applies to 'atomic' sentences. But do not modern formal systems with their full inductive definitions of formulas have greater complexity and precision than can be reconstructed from Aristotle's system embodying an only seemingly rudimentary account of complex sentences ? Against this worry we reply that the fact of the matter is that rigorous structural inductive definitions of sentences are not an exclusive possession of modern logic. There is strong evidence for such inductive definitions in Stoic Logic\cite{bob1} and specially  in  Book I 1173D-1176D of Boethius' \emph{De topicis differentiis} which will be discussed in detail in Section 9.  Finally we mention that perfectly 'modern' definition of the formulas of what we would now consider an extension of modal propositional logic (a multi-modal propositional logic)  is found in Ockham's \emph{Summa Logicae}\cite{ock}[79-80]. 
 
 \subsection{The Logic of the Prior Analytics}
 
 The four basic kinds of \emph{protasis} dealt with in the \emph{Analytics} can be seen as corresponding to four $(1,1)$-predicate term symbols. Syllogisms correspond very naturally to complex sentences in ESOL of a certain structure:
 
 \[ \forall A,B,C \phi(A,B) \& \psi(B,C) \rightarrow \rho(A,C)  \tag{Syll}\]
 
 where $\phi,\psi$ and $\rho$ range among such $(1,1)$- predicate term symbols (the order of the free variables is not necessarily the order of the free variables in their application to the corresponding symbol). For example, the Barbara syllogism corresponds to:
 \[ \forall A,B,C(\kappa AB \& \kappa BC \rightarrow \kappa AC)\]
 
 where  $\kappa XY$ is to be read as 'all $X$ is $Y$'.
 
From a formal point of view this is exactly what Aristotle is expressing and what Aristotle makes use of in the structure of the proofs in the \emph{Analytics}.  Aristotle certainly seems to assume the following reductive axiom (or definition) in the proofs:
 
 \[ \forall A,B(\kappa AB \leftrightarrow \forall x(Ax \rightarrow Bx)) \tag{A}\]
 
 Thus $\kappa$ is extensional. But Axiom (A) has to be tweaked as there is no direct evidence that Aristotle accepted $\kappa A A$.  To remedy this we introduce an extensional second-order equality predicate
 
 \[ A =_e B \leftrightarrow \forall x Ax \leftrightarrow Bx\]
 
 and then replace $(A)$ with

 \[ \forall A,B(\kappa AB \leftrightarrow \forall x(Ax \rightarrow Bx)\& A\neq_e B) \tag{A'}\]

As is well known, Aristotle seems to assume that all terms in the syllogistic have non-empty extension. But this evidently refers not to all terms (such as 'unicorn') but to those considered fundamental or relevant to scientific discourse (though given a term, the examination of whether this term has instances in patently a part of Aristotle's scientific methodology). Thanks to ESOL we can state this directly (using the $(1)$-predicate term symbol $\rho$) as an axiom
 
 \[ \forall A (\rho A \rightarrow  \exists x Ax) \tag{NEmp}\]
 
 In this Section we implicitly assume that quantifiers over monadic predicate variables  $X$ are subject to the restriction of $\rho X$.
 
 Similarly to $\kappa$  we can define $\mu XY$ to express 'some $X$ is $Y$'.  We have the following axiom or definition
 
 \[ \forall A, B(\mu AB \leftrightarrow \exists x(Ax \& Bx)) \tag{I}\]
 
 which states that $\mu$, like $\kappa$, is extensional. Again we mention the problem that we have to exclude $\mu A A$ and this leads us to variant
  \[ \forall A, B(\mu AB \leftrightarrow \exists x(Ax \& Bx) \& A\neq_e B) \tag{I'}\]
 
  We can follow the approach of \cite{pro} and introduce a second-order function to represent term-negation (\emph{de re} negation) which we denote by $A^c$. But thanks to our ESOL framework we can now define it directly
 
 \[ A^c = \lambda x \sim Ax \tag{c}\]
 
  We can define $\hat{\kappa} XY \leftrightarrow \kappa XY^c$ and $\hat{\mu} XY \leftrightarrow \mu XY^c$. Alternatively we could introduce $\hat{\kappa}$ and $\hat{\mu}$ directly with the expected axioms, for instance
 
 \[ \forall A, B (\hat{\mu} A B \leftrightarrow \exists x( Ax \& \sim Bx) ) \tag{O}\]

It is then straightfotward to derive all the assertoric syllogisms consisdered valid by Aristotle as well as to adapt Aristotle's own counterexamples to demonstrate the invalidity of those he considered invalid.

We now consider briefly how the modal syllogistic could be treated in the framework of ESOL. To this end we consider an additional second-order monadic function $\square X$. This takes a monadic predicate term $T$ and yields a monadic predicate term $\square T$. We posit the following axioms:

\[ \kappa \square A A \tag{1}\]
\[ \kappa \square A \square \square A \tag{2}\]
and the rule (K): from $\kappa A B$ we can derive $\kappa \square A\square B$. We also introduce the definition

\[ \sqcup A = \lambda x(\square Ax \vee \square A^cx)\]

and call the resulting system ESOL$^\square$. It seems very plausible that AML(see \cite{pro}) can be embedded into ESOL$^\square$ and all the corresponding theorems  be derived.

\section{The Logic of the Topics}

\subsection{Basic framework}

Aristotle's \emph{Topics} is generally considered to be an early and composite work, reflecting perhaps rather than specifically Aristotelian ideas, common dialectical procedures employed at the Platonic Academy. The purpose of the work is to serve as a debating manual. It is the part of the \emph{Organon} that has received less attention from scholars. But such neglect is undeserved as  the \emph{Topics} is quite a rich work and of considerable interest in its own right and it is impossible to do justice to all its aspects in the scope of the present paper. The \emph{Topics} transcends the strictly 'logical' concerns of the \emph{Analytics} : grammatical, semantic, ontological, pragmatic and epistemic aspects of debate are treated in equal measure.  The concern with ambiguity and synonymity amounts almost to an obsession, a fact that is likely explained by  the preoccupation with facing the rhetorical devices of the sophists. In this Section we confine ourselves primarily to the purely logical content of the \emph{Topics}.  
The \emph{Topics} is a debate manual, but for a very restricted and formalized kind of debate, concerned with a set of basic relations among terms\footnote{Unfortunately Aristotle's view of grammar and his theory of terms and concepts is only documented in a few scattered passages (but see chapter 4 of book IV which deals with relative terms and their governed cases).  A basic linguistic entity is the sentence or protasis.  A sentence is made up of two linguistic expressions connected by a certain propositional relation. The two expressions are called subject and predicate depending on their position in the propositional relation.  Only noun phrases can be a subject but adjectival and verbal constructions can be predicates. These can become subjects after a certain syntactic transformation. The original concept of 'verb' appears to have been that of predicate. A 'term' oros is what can be a subject.  However there is also a sense in which a verb is a term to which is added a temporal aspect.
The Stoic concept of kategorema seems to fit more closely to the concept of predicate. A predicate can be considered in itself, as unsaturated by any concrete subject. In Plato's \emph{Sophist} there is a discussion concerning 'sentences', 'nouns' and 'verbs', sentences being combinations of nouns and verbs, and all sentences must have a 'subject'. Plato gives the example 'Theaetetus, with whom I am now speaking, is flying'. This illustrates how Plato's concept of 'noun' corresponds to our 'noun phrase' and very likely to the Aristotelic term or 'oros'. There is very illuminating material regarding terms and the logical syntax of propositions ot be found in Boethius which is discussed in detail in Section 9.} (which the scholastics called \emph{predicabilia}), the most important being the species-genus relation, the property relation, the relation of difference (which combined with the species-genus relation yields the further relation of definition) and the relation of accident\footnote{It is not clear if this exhausts all the forms of predication. Indeed Aristotle often mentions predication without any further qualification and we are left to wonder if it assumed that this predication must fall within one of the classes above. Also 'accident' is understood in various different ways. But see the additional comments at the end of this Section.}. The second-order predicates $\kappa$ and $\mu$ of the \emph{Analytics} are a simplification and abstraction of these relations. Following \cite{prim} and the testimony of Boethius regarding a topic as a  'propositio maxima' we argue that a \emph{topic} is a sentence, generally in universally quantified form, expressing properties of these fundamental relations which can be used either to construct or refute a thesis involving one of the \emph{predicabilia}.
We argue that second-order logic is required to express the topics as well as the logical equivalences and inferences implicitly assumed to hold amidst related groups of topics\footnote{In 100a Aristotle gives the definition of 'syllogism'. The genus is 'logos'.  The difference is that in this logos certain premises being posited something different (heteron) from them  occurs (sumbainei) out of necessity. We might think that a syllogism thus corresponds to our modern notion of logical 'rule' (as in natural deduction or sequent calculi). The problem is that we would then have to reject rules such as $\Gamma, A \vdash A$ for $A$ is not different from the premises. Our position here is rather that 'syllogism' in this context corresponds not to the modern notion of 'rule' but to a universally quantified second-order sentence as explained further ahead. }. In our framework of ESOL  the \emph{predicabilia} are again $(1,1)$ predicate term symbols which we denote by $\gamma$ (species-genus), $\iota$ (property), $\delta$ (difference) and $\sigma$(accident).  But now they are no longer extensional, contrary to the case for $\kappa$ and $\mu$, that is, our relations cannot be reduced to 'first-order' expressions. However extension plays an important role in many topics: we need to express that the extension of one predicate is strictly contained in another or that they are equal. For example if we have $\gamma AB$ then the extension of $A$ is strictly contained in that of $B$, that is, $\forall x Ax \rightarrow Bx$ and $\exists x Bx\& \sim Ax$, but this last condition is not sufficient (contrary to the case of $\kappa$) to guarantee that $\gamma AB$.
 Also noteworthy is that the \emph{Topics} is not restricted to terms corresponding to  monadic predicate term variables, there are 'relation' terms as well, corresponding to predicate term variables of higher arity.  And we can wonder what for Aristotle the 'extension' of such as relation term would be. 
 
 The \emph{Topics} (which has traces of being a compilation of diverse, not always consistent, texts written at different times) consists of eight books. Book IV is dedicated to $\gamma$, book V to $\iota$ and book VI to definition (including a treatment of $\delta$). Important material relating to $\sigma$ is scattered throughout the first two books. Book IV is the most approachable from a formal  point of view and seems to also provide a transition to the \emph{Analytics}. For convenience we will in general write the topics as open formulas where the monadic predicate term variables are assumed to be universally quantified. The most common form of a topic is
 
 \[ \forall A, B(  \phi(A,B) \rightarrow \bigcirc A B)\]
 
 where $\bigcirc$ represents a predicabilia or the definition relation $h$ (defined further ahead) or else the negation of a predicabilia  or the negation of the definition relation.  $\phi(A,B)$ is any formula with free predicate term variables  $A$ and $B$. The formula $\phi(A,B)$ will in general contain second-order quantifiers as well as first-order quantifiers. It follows that these topics can also be written
 \[ \forall A, B( \bigcirc A B \rightarrow  \phi(A,B))\]
 
 For the genus-species relation $\gamma$ there is evidence that Aristotle posits:
 
 \[\gamma AB \rightarrow \forall x(Ax \rightarrow Bx) \& \exists y(By \& \sim Ay)\quad \text{($\Gamma$)}\]

 We have the difference relation $\delta AB$ read as '$A$ is a difference of genus $B$'\footnote{ A careful reading of the text (specially Book VI) reveals that Aristotle uses 'difference' in three different ways: the  difference of a species (relative to a genus or a definition), coordinate differences (of the same genus) and differences that 'belong' to a given genus.  We consider the third to be the most fundamental and define the two previous notions of difference in terms of this last one.}.
 For property\footnote{The topics of property set the stage for the topics of definition which is the culmination of the \emph{Topics}. The notions deployed for property are basically the same one as those deployed for definition but in a more straightforward way. There is a much greater laxness in what can count as a property as compared to definition.  We might say that property represents of a more archaic form of definition before the requirement of the a regimentation into a (proximate) genus and difference was adopted. 
 	The basic idea is that the property of a term $A$ cannot involve $A$ itself (non-circularity), cannot employ less known terms  and cannot contain
 	redundancies. } $\iota AB$, read '$B$ is a property of $A$' , we have
 
 \[\iota AB \rightarrow \forall x(Ax \leftrightarrow Bx)\quad \text{(I)}\]
 
 None of these implications are equivalences\footnote{There is no doubt that Aristotle conceives of co-extensive terms which are intensionally different (that is, they are essentially different terms, not only syntactically different). Likewise the terms 'One' and 'Being' are described as being predicated of all other terms and yet intensionally distinct.}.
 
 We could also venture the following axiom for accident\footnote{Aristotle distinguishes between separable and inseparable (or essential) accident. But this seems to make sense primarily for predication of an individual.} $\sigma AB$ read '$B$ is an accident of $A$':
 
 \[\sigma AB \rightarrow \exists x,y(Ax \& Bx \& Ay\& \sim By)\quad \text{($\Sigma$)}\]
 
 However the version of accident $\sigma'$ in a passage in 120b suggests we need a modal version of ESOL. An accident is simply a predicate which can both be and not be predicated of an individual which we could also read in ESOL$^\square$ as
 
 \[ \sigma' A \leftrightarrow \exists x( \lambda x(\lozenge A \& \lozenge  A^c)) x \]

 Consider axiom ($\Gamma$). This is an example of a topic. It is neither a topic of refutation nor directly a topic of construction - but an intermediate kind, a make sure that so-and-so topic. But from ($\Gamma$) Aristotle infers a powerful and convenient topic of refutation:

 \[ \forall A,B(\exists x(Ax\& \sim Bx) \rightarrow \sim \gamma AB) \]
 
 which is topic 5 of Chapter 1 listed in appendix 1.
 
 Although the topics are usually listed in an \emph{ad hoc} fashion it is clear that they are organized into closely knit units bound by logical equivalence or implication. ESOL s required both to formulate the topics and to explain how Aristotle might have established the inferences and equivalence therein. We might call the working out of these implicit inferences the 'metalogic' of the \emph{Topics}.
 Our objective here is not to undertake an exhaustive formal interpretation of all the logical topics in   ESOL. Rather it is to go through a selection of key topics mainly from book IV but also from books V and VI which provide evidence for the (extended) second-order nature of the metalogic.  For convenience instead of referring directly to the corresponding passage in the text of the \emph{Topics} for each topic discussed we make use of the lists compiled by Owen in \cite{ow}: a formalized version of a portion of them is found in Appendix 1. 
 Let us start with book IV and $\gamma$. We find therein topics which state directly or are stated to be inferred from axioms which express the fact of $\gamma$ having  transitive, symmetric and anti-reflexive properties:
  \[ \gamma AB \& \gamma BC \rightarrow \gamma AC \]
 \[ \sim \gamma AA\]
 \[ \gamma AB \rightarrow \sim \gamma BA\]
 
 as can be verified in appendix 1.
 The first topic is clearly an ancestor of the Barbara syllogism of the \emph{Analytics}. Perhaps $\mu$  descends  from $\delta$ since the same difference can operate on distinct genera so that $\delta AB$, $A$ is a difference of genus $B$ means that 'some $B$ are $A$' (a topic discussed further ahead has certainly some connection to the syllogism Darii).  There is the important tree-like property:

 \[  \gamma AB\& \gamma AC \rightarrow (\gamma BC \vee \gamma C B) \]
 
 which is the first topic in chapter 2 of book IV. Aristotle refers to $\gamma AB$ also as '$A$ participating in $B$' or '$A$ receiving  the definition of $B$'. Aristotle mentions important axioms in the context of the derivation of refutation or establishing topics. This gives us insight into the metalogic of the topics. For instance in 122b Aristotle infers a topic using the tree-like property and the topic can itself be represented as an inference in sequent form:
 
 \[  \forall AB (\gamma AB\& \gamma AC \rightarrow (\gamma BC \vee \gamma C B)), \gamma DE, \gamma DF, \neg \gamma EF \vdash \gamma DE \]
 
 which of course makes use of the 'transitive' property of $\gamma$.
 
 Aristotle has the concept of \emph{proximate genus} which can be easily expressed in ESOL
 
 \[  \gamma^\star AB \leftrightarrow \gamma AB \&\sim \exists C(\gamma AC \& \gamma CB)\]
 
and so too can the important postulate that every genus must have more than one proximate species (a species of which it is the proximate genus):
 
 \[ \gamma^\star AB \rightarrow \exists C(\gamma^\star CB \& C\neq A)\]

But the problem here is to account for the 'equality' relation which evidently is not the same as extensional equality $=_e$ but is stronger:

 \[ A = B \rightarrow (\forall x Ax \leftrightarrow Bx) \tag{$I=$}\]
 
 but we do not have (as we argued previously) the reciprocal implication. We could tentatively posit that we have a second-order version of Leibniz's law for $=$(see the topics in 152a which suggests particular cases of this law). In particular in book VII 152a31-33 we read (tr. W. A. Pickard):\\
 
 \emph{Again, look and see if, supposing the one to be the same as something, the other also is the same as it: for if they be not both the same as the same thing, clearly neither are they the same as one another.}\\
 
 Aristotle appears to be stating that $A = B \rightarrow (A = C \rightarrow B = C)$ since \[\sim ( A = C \rightarrow B = C) \rightarrow \sim A = B\]

 Axioms ($\Gamma$) and (I)  are stated as a topics.  So too are often the counter-positives of some of the topics considered above 
 
 All the relations above are  exclusive of each other. For instance we have the topic $\sigma A B \rightarrow \sim \gamma AB$ - but this can be derived from ($\Sigma$) and ($\Gamma$). We also have that $\sim \Sigma AA$. There is a topic which seems to state that if $\gamma AB$ then $B$ cannot be the difference of a genus of $A$:
 
 \[ \gamma AB \rightarrow \sim \exists C(\delta BC \& \gamma AC)\]

 If $\gamma AB$ then there must be some difference $C$ of $B$, $\delta CB$ which 'is said' of $A$. Aristotle has a kind of unspecific general kind of 'predication' wherein it is not obvious whether it is supposed to fall under one of the predicabilia or not. We take this unspecified predication to be just extensional containment $A\subset B \leftrightarrow \forall x (Ax \rightarrow Bx)$. 
  We might express this previous  topic as
 \[ \gamma AB \rightarrow \exists C(\delta CB \&  A \subset C)\]
 
 From the above we derive the refutation topic:

 \[ \forall C(\delta CB \rightarrow \sim A\not\subset C) \rightarrow \sim \gamma AB\]
 
 which is listed as number 11 of Chapter 2 in appendix 1.
 We find in addition the following interesting topic:
 
 \[ \gamma AB \rightarrow \forall C (\delta CA \rightarrow \delta CB) \]
 
 If bipedal is a difference of mammal it is also a difference of animal.

 We will not analyze here the question of supreme genera and infima species but the order-theoretic properties involving $\gamma$ could be easily expressed in ESOL though other implied properties such as the finiteness of ascending chains of subordinate genera would have to make use of an extension of ESOL.

  We have not addressed the ontological categories which feature in the \emph{Topics}. The order on terms induced by $\gamma$  can in general  have more than one maximal element. And we can ask if it is correct to interpret first-order individuals as individual substances. There is the problem that Aristotle seems to distinguish different kinds of situation for which $Ax$, for which a predicate term $A$ holds of an individual substance $x$, for example, if $A$ is predicated according to essence or as an accident. This would suggest that we need  an extra predicate term symbol with  second-order and first-order arguments $\gamma^i Ax$ (for which $\gamma^i Ax \rightarrow Ax)$ to make such ontological predication distinctions . A general approach in ESOL to the theory of predication in the \emph{The Categories} might be something as follows (we will not discuss here the question of whether some of the other categories besides relation are in fact also relations: knowledge is a state of a subject but also, it seems, a relation). We postulate nine different categories of predication $\mathcal{P}^iTx$ where we put $\mathcal{P}^1 Tx \leftrightarrow Tx$ as corresponding to \emph{ousia} while the others correspond to the other Aristotelian categories. The tenth will correspond to relation and must split into signatures of the form $\mathcal{P}^{10}R^2xy$ and  $\mathcal{P}'^{10}R^2xT$. For consider 'knowledge' and 'knowledge of colours'.  Knowledge corresponds to the state of a subject but also to a relation, and this relation can regard not one or many primary substances but other categories: if a subject $a$ knows the colour 'blue' then this is not the same as knowing all blue substances. For another perspective on the categories see the additional notes in the final subsection of this Section.

 A central role in the \emph{Topics} is played be a kind of 'semantic' term negation operator (opposition operator),  operating within the proximate species of a given genus. This must be a second-order function applied (primarily) to monadic predicate terms. Its properties are
 
 \[ \gamma^\star AB \rightarrow \gamma^\star A^\bullet B^\bullet\]
 \[ A = B \rightarrow A^\bullet = B^\bullet\]
 
 and we also assume $A^{\bullet\bullet} = A$.
 
 For instance the opposite of the color white is black. A proximate species of a genus may not have an opposite; this situation is expressed as $A^\bullet = A$. Noteworthy topics are 
 \[\gamma AB \rightarrow \gamma A^\bullet B^\bullet\]
 and also
 \[ B = B^\bullet \& \gamma AB \rightarrow \gamma A^\bullet B\]
 
 as can be verified in appendix 1.
 
 In Book V we encounter the topic
 
 \[ \iota AB \rightarrow \iota A^\bullet B^\bullet \]

 In order to formalize the theory of definition we cannot do without a  binary function on monadic predicate terms which corresponds grammatically to a modifier in a noun-phrase, for instance an adjective attached to a head noun. We study in more detail the grammatical aspect of definition terms in the next subsection.  We denote such a function by means of concatenation and impose the axiom which is implicitly assumed by Aristotle:
 \[ \forall x((AB)x \leftrightarrow Ax\& Bx) \tag{Concat}\]
 
 Since we do not want semantically meaningless constructs like 'green numbers' it makes sense to define a further second-order monadic predicate allowing us to express non non-empty extension (cf. also (NEmp) in Section 3):
 
 \[ AB\downarrow \leftrightarrow \exists x(AB)x\]
 
 so that
 
 \[ \delta BC \rightarrow BC \downarrow\]
 
 The fundamental property of concatenation is that it lets a difference carve out a species:
 
 \[ \delta D A \rightarrow \gamma (DA) A \]
 
 We can define definition\footnote{Aristotle also uses the term 'logos' for definition (cf. 141b3).} according to Book VI as
 
 \[h AD \leftrightarrow \exists B, C( D= BC \& \delta BC \& \gamma^\star AC)  \]
 
 We deduce immediately from $(I=)$ that $h A D \rightarrow (\forall x Ax \leftrightarrow Dx)$. 
 We may now formulate a fundamental axiom of Aristotelian logic, that every concept (in this case monadic predicate) has a definition:
 
 \[ \forall A \exists D \,hAD \tag{h}\]
 
 From this and the axiom for our concatenation operation (Concat) we can deduce the that immediate species of a genus have disjoint extensions:
 
 \[ \forall AB\, \gamma^\ast AG\& \gamma^\ast BG  \rightarrow A = B \vee (\sim \exists x(Ax\& Bx)) \]

 It is important to note that equality here is not to be read syntactically but semantically. Thus $D = BC$ is to be read as 'the concept $D$ has as parts concepts $B$ and $C$', the paradigm being 'rational animal'.  It is a fundamental problem whether a term can be said to be equal to its definition (with this notion of equality) and whether Leibniz's law holds for this notion of equality. It could be that two distinct kinds of equality are required. It is true however that  definition though the substitution of a term by its definition features in several topics\footnote{152a-152b is very relevant for term equality. Aristotle announces topics which express the transitivity of $=$ as well as topics recalling Leibniz's law. In 153a he seems to say that there is equality between a term and its definition but this is in itself not enough to show that it is a definition. The \emph{definiens} must satisfy further conditions. 
 	The passage 140b32-141a1 is also interesting.
 	Aristotle states that a definition is identical to the term it defines:
 	\emph{esti de tauton tô(i) anthrôpô(i) zô(i)on pezon dipoun}. However in 143a15-28 where Aristotle appears to admit that either giving the immediate genus and the difference or a non-immediate genus and all corresponding differences (the 'intension') as both valid forms of definition which leads us to inquire exactly what the nature of the purported 'uniqueness' of the definition.  Does this means that both these forms are to be regarded as the same definition ?}.
 From this last definition we can evidently derive (in ESOL) numerous refutation topics - and these are indeed topics in book VI. As Aristotle says: it is much easier to refute a definition than to establish one.
 
 A standard view would be that if $hA(DC)$ then $C$ must be uniquely determined (this in fact follows from the above axioms). There is however a problem suggested by the passage 148b dealing with 'complex terms'\footnote{this notion of 'complex term' appears to be primarily grammatical and not logical or semantic as is our concatenation operation, however the precise relation between the two remains illusive as does the interpretation of the passage in question. Are we to understand the crucial non-circularity property of the definiens in semantic or syntactic terms ? But semantically does it even make sense to have a situation in which $A = AB$ ? }.
 
 Up to now 'difference' was defined relative to a genus, that which carves out the species of a genus. We now wish to define what it means for a difference to be the difference of its species:
 
 \[ \Delta A D \leftrightarrow \exists C(\gamma^\star A C\& \delta D C \& h A (DC)) \]
 It is  interesting how definitions relates to the opposition operator. If we have a definition $DG$ of $A$ what can we say in general about the definition of $A^\bullet$ ? In 153b Aristotle says that there are three cases to consider for the definition for $A^\bullet$ depending on whether $G^\bullet = G$ or $D^\bullet = D$. The possibilities are either $D^\bullet G^\bullet$, $D^\bullet G$ or $D G^\bullet$. For example justice is virtue of the soul and its opposite injustice is the vice of the  soul, soul having no opposite.

 \subsection{On relations}
 
 We now discuss Aristotle's relation terms concerning which Section 4 of book IV of the \emph{Topics} is of particular interest. A natural approach is to represent  relation terms by binary predicate terms in ESOL. But then we need a corresponding version of $\gamma$ for such arguments, there is no way to deal uniformly with the predicabilia when predicate terms of different arities are involved. This is also enforced by the topic in book IV Section 1 (120b) which requires that species and genus belong to the same category, although Aristotle admits later in Section 4 that there might be exceptions (124b) as for the case of 'grammatical' (monadic) as a species of the genus 'knowledge'(relation). In our case the application of the aforementioned topic (stated also explicitly in (124b)) is simply that there must be two species-genus relations, one $\gamma$ for monadic predicates and another $\gamma_2$  for binary predicate terms. We need to postulate analogous axioms analogous to those for $\gamma$ such as 
 \[ \gamma_2RS \&  \gamma_2 ST \rightarrow \gamma_2RT\]
 
 where $R,S$ and $T$ are binary predicate term variables. The logic of relation terms is quite complex in Aristotle and bound up with Greek syntax in a way which renders interpretation difficult. A relation term $R$ has its 'reciprocal' or 'inverse' which we can interpreted by a  second-order function  $R^\ast$. It is a challenge is to interpret the passage in (124b) about the double ($D$) being the double of the half.
 It is plausible that we are in this case in the presence of relations for which
 
 \[ \forall x,x',y,y' (Rxy \& Rxy' \rightarrow y = y') \& (Rxy\& Rx'y \rightarrow x = x') \]
 
 One way of interpreting 'the double is double the half' would be
 
 \[  \forall x Dx (\iota y D^{\ast}xy)\]
 
 so we can read the topic in (124b) as 
 
 \[  \gamma_2RS \rightarrow  \forall x Sx (\iota y R^{\ast}xy) \]
 or in its refutation version
 \[  \sim\forall x Sx (\iota y R^{\ast}xy) \rightarrow  \sim \gamma_2RS \]

 In the passages that follow Aristotle discusses the necessary agreement between the case  governed for terms in the $\gamma_2$ relation. For example if $\gamma_2RS$ then if $R$ takes the genitive then so does $S$. More interesting from a logical point of view is the passage mentioning the number of complements, this suggests that Aristotle was conscious of the necessity of the general agreement between the arity of the predicate terms in the species-genus relation, so that we should consider a whole sequence of such relations $\gamma_1 = \gamma, \gamma_2,...,\gamma_i,....$  where $\gamma_n$ serves for $n$-ary predicate variables.  Aristotle uses the example of the trinary relations 'donation' and 'debt'.
 If we look at the early and middle Platonic dialogues we find that relation plays a much larger and sophisticated role than in does in Aristotle's \emph{Organon}. The logic of relations does not appear to be very worked at our and was perhaps consciously and deliberately neglected. Indeed there are some natural questions that could be asked (both now and in antiquity). For example: how do definitions work for relations ? What are the differences and how are they joined to the respective proximate genera ?
 The search for the definition of relatives  figures in the Platonic dialogues, for example in the \emph{Lysis} the definition of friendship is sought after.  A definition of a relation can necessarily involve sentences of arbitrary complexity involving other relations which do not fall under the general scheme for monadic predicate terms of book VI of the \emph{Topics} (at least for its most simple familiar cases). Indeed traces of such unorthodox definitions are to be founded scattered throughout this work\footnote{for example the definition of Justice as Courage and Temperance (150a-151b).}. Even the most basic relations such as that of 'being a paternal grandfather'   do not seem to be organized neatly into a proximate genus $+$ difference format. In the \emph{Lysis} the binary relations of friendship $\Phi$ and love $L$ are considered. For example the following hypothetical definition is put forward:
 \[ \Phi x y \leftrightarrow Lxy\& Lyx \]

 The \emph{Lysis} deploys relatively complex statements and reasoning about relations - and none of the proposed definitions seem to fit the genus $+$ difference template. We can wonder why Aristotle did not take this more into consideration in the \emph{Topics} and \emph{Analytics}. However there are certain interesting passages which shed light on Aristotle's view of relations. Consider the passage: 145a14-15: \emph{tôn gar pros ti kai hai diaphorai pros ti}. Aristotle claims that  differences of the relation 'science' such as 'practical', 'theoretical' and 'productive' are themselves relations. It is not immediately clear how to express this in ESOL (but see the discussion further ahead). The discussion in 142a is interesting for it appears to concede that the rule of defining the less know by the better known may fail for relatives such as 'half' and 'double'. Such relatives must be defined in terms of their correlatives, the double in terms of the half  (\emph{to diplasion hêmiseos}, 142a26).

 As already mentioned, in Appendix 1 we list some of the most important 'logical' topics in Book IV including (but not limited to) those discussed in this Section. We now investigate the structure of the noun phrases used as definitions, which is relevant also to a refinement of our 'concatenation' operation.
 
 \subsection{The grammatical structure of definition terms}
 
  Classical Greek grammars do not give us much in the way of an account of the structure of noun phrases.  Thus we will use the concepts employed in modern linguistic treatments of English grammar such as \cite{leech}.  A noun phrase is usually considered as having the structure: determiner + premodifier + noun + postmodifer and relative clauses are included as a type of postmodifier. Obviously the case of being a pre- or postmodifier is not the same in ancient (Attic) Greek and modern English but we do not consider the position of a modifier expression relative to the head of the noun phrase to be of importance for our purpose, only the fact that it is a modifier.  
 Definitions in the \emph{Topics} are noun phrases and as such we must investigate their structure. Questions we could investigate in the future include: how are we to understand grammatically the genus+difference template endorsed by Aristotle ? Does the genus always correspond to the head noun in the structure of the noun-phrase ? Does Aristotle allow definitions that do not fit this regimentation ? How do definitions of qualities differ grammatically from those of secondary substances ? What is the grammatical structure of the definitions of relatives ?  For a term of a given category what are the constraints on the category of the difference ? And finally: is ESOL adequate to capture the logical content of all the grammatical forms of definitions considered by Aristotle ?
 In order to examine some of these questions it would be useful to list all the definitions (both accepted and rejected) discussed by Aristotle in book VI of the Topics but we list here only those definitions appearing in chapters 1 -5. We use the edition of the Greek Text from \cite{topiques2}. We present the definition noun phrase and then in parenthesis the noun of which the phrase is proposed of as a definition\footnote{Aristotle often presents both the term and its definition in the accusative case following the construction: X (acc.) Y(acc.) if (or is) defined. We have endeavored in most cases to revert the term and the definition back to the nominative case.} (whether accepted or not by Aristotle: most of the examples are rejected but for reasons which do not involve their grammatical structure).
 \emph{
 \begin{enumerate}
 	\item 139b20: agôgê eis ousian (genesis)
 	\item 139b20-21: summetria thermôn kai psukhôn (hugieia) 
 	\item 139b32-34 : ametaptôton (epistemê), tithênê (gê), sumphonia (sôphrosunê)
 	\item 140a3-5:  ophruoskios (ophthalmos), sêpsidakes (phalaggion), osteogenes (muelos)
 	\item 140b1-5: arithmos autos hauton kinoun (psukhê)
 	\item 140b7-8: hugron prôton apo trophês apepton (flegma)
 	\item 140a35 zôon pezon dipoun tetrapêkhun to epistêmês dektikon (anthrôpos)
 	\item 140b23-24: zôon pezon dipoun tetrapêkhun (anthrôpos)
 	\item 140b27-28 oreksin hêdeos (epithumia)
 	\item 141a7 horistikên kai theôrêtikên tôn ontôn (phronêsis)
 	\item 141a10  sterêsis tou kata phusin thermou (katapsuksis)
 	\item 141a16  elattôsis tôn sumpherontôn kai dikaiôn (epieikeia)
 	\item 141a19-20 epistêmê tôn hugieinôn zô(i)ôn kai anthrôpôn (iatrikê)
 	\item 141a19-20 eikona tôn phusei kalôn kai dikaiôn (nomos)
 	\item 141b21 epipedou peras (grammês), stereou peras (epipedos)
 	\item 142b1 astron hêmerophanes (hêlios)
 	\item 142b3-4 hêliou phora huper gês (hêmera)
 	\item 142b8 to monadi meizon artiou (peritton)
 	\item 142b11-12 arithmos ton dikha diairoumenon (artion), heksin aretês (to agathos)
 	\item 142b24-26 to ekhon treis diastaseis (sôma), to epistamenon arithmein (anthrôpos)
 	\item 142b30-31 to epistêmên tou grapsai to hupagoreuthen kai to anagnônai (grammatikê)
 	\item 143a3 tou noson kai hugieian poêsai (iatrikê)
 	\item 143a16-17 heksin isotêtos poêtikên ê dianemêtikên tou isou (dikaiosunê)
 \end{enumerate}}

The above list provides  evidence that a great variety of noun phrase constructions (and in particular with regards to the structure of the propositional phrases appearing as postmodifiers) were allowed by Aristotle. 
 In (1) and (2) the genera seem to be \emph{agôgê} and \emph{summetria} while \emph{eis ousian} and \emph{thermôn kai psukhôn} would correspond to prepositional noun phrase postmodifiers in English. In (6) it appears that the head or genus is \emph{hugron} which is affected by two modifiers, the adjective \emph{apepton} and the non-finite clause equivalent to a relative clause \emph{prôton apo trophês} were the infinitive (probably 'generating') is omitted.
 While a common noun + adjective definition for secondary substances can be construed as a conjunction of monadic predicate constants $G$ and $D$ yielding a definition term $\lambda x(Gx\& Dx)$ the situation is more delicate for other kinds of postmodifiers. One approach is to consider the case in which the genus or head is  a relation. For instance \emph{agôgê} (coming to be, arising) is modified to 'coming to be according to essence'. If $A$ represents the relation of coming to be according to something than the last definition could be represented by $\lambda x Ax\top$ where $\top$ is the supreme genera for \emph{ousia}.
 
 We cannot avoid mereology in (2): given a substance $a$ we need to form the sub-substances of the cold things in it $c$ and the warm things in it $w$. Could this be expressed by a term-forming operator syntactically analogous to set comprehension: $c = \{ x \epsilon a : Cx\}$ and $w = \{x \epsilon a: Wx\}$ where we also have the analogue of $\cup$ ? If \emph{summetria} is a relation $S$ then we interpret $Sxyz$ as meaning that $x$ has symmetry of $y$ and $z$.  Thus our definition of \emph{hugieia} (health) would look like
 
 \[  H = \lambda y(Sy \{ x \epsilon a : Cy\} \{ x \epsilon y : Wx\})\]
 
 In (20) 'that which has three dimensions' we have an example of a noun-phrase having as head a demonstrative pronoun with a relative clause postmodifier.
The passage 140b32-141a1 is interesting. Aristotle's point here seems to be that postmodifiers are not predicative so that a term may be predicated of a noun phrase (\emph{kata tou sumpantos}) with a modifier employing the same term without redundancy (\emph{hapaks hê kategoria ginetai}), for example: a biped animal is biped\footnote{This might be interpreted as Aristotle pointing out the logical-grammatical difference between predication of a concept to a subject and the concept used to specify the subject. }.
Finally we should not fail to mention that ESOL is adequate to treat propositional attitudes through the use of second-order predicates involving $0$-ary predicate terms (which correspond to propositions) and the $\lambda$-construction for an empty list of variables as standard in Second-Order Logic. For example $a$ knows that $Ab$ could be expressed as $Ka(\lambda(Ab))$. It would be interesting to look for instances of propositional attitudes in the \emph{Organon}.

\subsection{Additional remarks on Book I of the Topics}

Chapter 1: Aristotle's definition of 'syllogism'  here is quite general (and should not be confused with the 'syllogisms' of the Analytics) and a good translation would be 'inference', the kind of inference represented by a sequent in the sequent calculus (with the cut-rule) with one formula on the right $A_1,A_1,....,A_n \vdash B$.  Aristotle's 'true' or 'primary'  things are a sequents of the form $\vdash A$.

Chapter 2. But what is dialectic in the Topics ? Does it investigate the axioms of the particular sciences themselves (which cannot be investigated in those sciences) ? The passage 101b1-3 is mysterious: dialectic is a critical path possessing  the 'beginning of all methods'.

Chapter 4.  The protasis. Each protasis indicates (is made up from) either property, genus or accident. Difference is strangely classed as 'pertaining to genus (generic)'. And   'problems' can be constructed from every protasis by changing the 'mode'.  Can we devise a term-forming operator corresponding to interrogation ?

Chapter 5. This is a very important Section.  We have a 'logos' which 'semainein' (signifies).  It is not clear if the logos is meant as a mere signifier or as the sign (signifier + signified).  Here logos is contrasted with onoma.  Apparently this corresponds to the difference between simple and complex terms. We see that there are definitions of complex terms.  An important question is: can we accept a definition consisting of a simple term ? Here Aristotle hesitates but admits that such a protasis are at least useful for definition. A fundamental concept is that of 'antikategorein' (to be convertible with).  The later Fregean distinctions between Sinn and Bedeutung  as well as between concept, object (and extension) are all present in the Topics.  A is convertible with B for Aristotle if A and B are predicates with the same extension (but not necessarily with the same meaning).  Aristotle's formula in 102b20-23 is : if (it) is A then (it) is B and if (it) is B then (it) is A.  Clearly definition and property have the same extension but different meaning (they both do not signify essence).  The rest of the discussion is valuable for elucidation of 'accident' and how it overlaps with relative (and temporary) property. The example of the 'only man sitting' (in a group) suggests a connection to definite descriptions of individuals.

How did Aristotle distinguish between: 'all men are animals',  'man is an animal', 'animal is the genus of man', between extensional and intensional predication ?  In the above chapter Aristotle seems to give the rule: from $\Gamma \vdash \forall_{x:A(x)} B$ we can derive $\Gamma \vdash \exists_{x:A(x)} B$.   In the beginning of chapter 6 of the second book we find the Stoic exclusive disjunction.  It is patent that 'connectives' for Aristotle are just as much term operators as operators on propositions.

Maybe  a dependent, alternative, version of existential quantifier introduction is required.  We could reserve universal quantification for explicitly 'distributional', extensional quantification (such as in the the expression of the topics itself) and in other case use $\gamma$, etc. And $\exists^\gamma_A B$ would a predicate meaning $\exists C\prec B \gamma AC$. Thus our new version of the existential quantifier rule (which is a topic at least implicit in book 2, 109b) would look like this: from $\gamma AC$ and $C\prec B$ deduce $\exists^\gamma_B A$.

Chapter 7. Peri Tautou (sameness, identity).  There is a kind of qualitative sameness considered here and the example of the water drawn from a given spring is noteworthy. What is Aristotle's 'arithmetical identity':  the name being many and the thing being one (103a9-10)?  Here Aristotle may be interpreted as postulating that identity is not a primitive notion but a polysemic and it to be defined in terms of either homonymity, definition or property  or even accident. This rules out any extensionality (Frege's Law V).

Chapter 8.  Definition consists of genus and differences and these are said to be 'in'  the definition. We still need to investigate further the precise relationship between difference and accident. 

Chapter 9.  For Aristotle 'kategoria' means 'predicate'.  What is the relationship between onoma, protasis, logos and kategoria?  Category in the ordinary sense is actually 'genus of predicate'.   The ten classes seem to be genera both of predicates and of things in general  - an ontology.   We can raise the question of the definition of objects that do not belong to the category of substance. Also about Aristotle's view on statements of the form 'A is A'.  Aristotle is stating here that if the thing and predicate belong to the same class then we have an essential predication, otherwise we do not.  But how can we accept Aristotle's example of predicating 'man' of a given man being a predication according to essence ?  How can 'white is white' signify essence? 

Chapter 10. This chapter offers us the rudiments of a new kind of intensional logic: a doxastic logic, and is of considerable interest.  Although we can think of ways of making sense of a doxastic de re modality, it is perhaps better to consider first a de dicto second-order relation $\xi AB$ to be interpreted as: $B$ is generally held to be predicated of $A$ - it is plausible that $B$ is predicatedd of $A$. 

Chapter 13.  Differences of meaning of a term and a term qua term can be objects themselves of propositions.  To formalize the Topics we thus may need a semantic identity relation.

Chapter 15.  It would be interesting to investigate formal systems in which each term is assumed to be interpreted as having a set of references (and meanings) rather than one (or none).  This is the kind of polysemic logic that looms large in the Topics. A kind of semantic set theory, perhaps.  The task is to construct expressions which are singletons and to detect them within the formal logical and grammatical rules of the system. Aristotle must accept that there is a notion of semantic identity (which is not the same as that of 'antikategorein' or extensional equivalence). 

\section{Formalizing book V}
In this Section we focus on book V and discuss some of the challenges facing the ESOL approach as well as its versatility in solving them.

It is probable that the elaborate theory of definition in book VI grew initially from the theory of property (idion) . Some textual evidence can be adduced in favor of this (cf. 154b). There are passages in which Aristotle relaxes all the complex structural requirements on the definiens term, even simple terms are allowed to be definitions; on the other hand a  passage in 132a indicates that property must have a genus + specifications structure just like definition. How are we to understand that 'terrestrial bipedial animal' indicates the essence of man but that (132a) 'tame (hêmeron) animal'  or 'animal gifted with knowledge' somehow does not ?  It might be tempting to view book V of the Topics, dealing with property, to be both among the oldest portions of the work and  even as having constituted an independent treatise.  It is true that definition, even difference as well as predication 'according to essence', are mentioned (and contrasted with property) in book V,  but not in a detailed way, rather the passages in question could easily be viewed as later interpolations or revisions. We propose that book V be studied on its own in a more-or-less self-contained way.

In the most basic sense, property is related to co-extensionality as in $(I)$. But this is not all. Originally property, idion,  seems to have been subordinated to  pragmatic, epistemic purposes. Property was a kind of  'distinguishing mark' (perhaps not unlike \emph{linga} in the Indian tradition)  which served to make an object known.  Thus it is natural to suppose that the property is somehow conceptually simpler than that which it is a property of,  as in certain marks distinguishing different species of plants. 

For Aristotle co-extensionality is a necessary but not sufficient condition for being a property (cf. also 153a in book VII: (extensional ?) equality is not enough for there to be definition). Thus being a property is intensional. If $\iota AB$ we cannot substitute co-extensional terms for $A$ or $B$ \emph{salva veritate}.  Indeed Aristotle's predicate of terms ' being better known'  (see 129b-130a) is also intensional. We denote this situation by using a second-order predicate tem symbol  $A\prec B$ to mean that the term $B$ is better known than $A$. $\prec$ is a transitive.  We write $A\sim B$ if $A$ and $B$ are simultaneously known (this can be defined as $A\nprec B\& B\nprec A$).  We meet the following rule in 131a: $ A \sim A^\bullet$.

There is a series of syntactic and semantic  conditions imposed on the property term. This theory involves implicit syntactic-grammatical notions which unfortunately are never made explicit in the surviving texts. Here are, in a modernized language, some of the conditions. For $B$ to be a property of $A$ it is required that:
\begin{enumerate}
	
	\item  the term $A$ does not occur in $B$ (non-circularity) (131a)
	
	\item  the term $B$ is different from the term $A$ (non-identity) (135a)
	
	\item the term $B$ is not a conjunction $\lambda x(B_1x \& B_2x)$ in which $B_1$ or $B_2$ are already properties of $A$ (irreducibility)
	
	\item $B$ is not a universal predicate or does not make use of universal predicates (non-redundancy)(130b)
	
	\item the term $B$ does not make use of terms lesser known than $A$ (epistemic priority)(129b)
	
	\item the term $B$ does not have repeated subterms (including via anaphora) (130b)
	
\end{enumerate}

Perhaps 3 can be replaced with\\

3'.the term $B$ does not strictly contain a term $B'$ which is a property of $A$\\

As a corollary of 5 we obtain that: $B$ cannot contain $A^\bullet$ (131a). 1 is actually strengthened by Aristotle to the requirement that $B$ does not even contain species of $A$.

Modern systems of formal logic do not usually possess a mechanism to talk about their own syntax (at least not directly). And the conditions above do not seem to be primarily semantic (as our operation $AB$) but syntactic. But there is nevertheless a certainly a semantic aspect involved and we could try to formalize the semantic aspect in ESOL by means of our 'concatenation' operation $(AB)$. For instance we could write

\[A < B \leftrightarrow \exists C(A = BC) \]

and $A\leq B$ for $A < B \vee A = B$.

and formulate weaker versions of some of the above topics as

\begin{enumerate}
	\item $\iota AB \rightarrow A \nless B$
	\item $\iota AB \rightarrow A \neq B$
	\item $\iota AB \rightarrow \sim\exists C(\iota AC  \& C < B)$
	\item $\iota AB \rightarrow \sim \exists C(C\leq B \& \forall x Cx)$
	\item $\iota AB \rightarrow A\prec B$
\end{enumerate}

It is difficult to formalize 'universal predicate' for what kind of predication (out of the predicabilia) is this referring to ? Aristotle does not accept that a universal predicate can be a genus thus this can't be predication understood as $\gamma$. The interpretation we take here is a weak version 4 referring to universal predicates according to extension. As our semantic concatenation operation arguably satisfies $AA = A$ it does not seem possible to express 6. 

None of these conditions would we rejected, for instance, by a modern mathematician who wishes to characterize (but not necessarily define) a class of objects. Even if the 'better known' relation is questionable, characterizations involving theories disproportionately more elaborate and distant than those of the the original class of objects will not always be welcomed.

Does it make sense to speak of a property of an individual ? How does one extrapolate from predication of individuals to a property of their species ? Aristotle introduces the construction: '$B$ is said of all individuals $a$ such that $Aa$ \emph{qua} $A$' (132b). How are we to interpret this ?  Aristotle in 132b speaks of the 'onoma' and the 'logos'.  But we could give a simple interpretation of the refutation topic in question in terms of extensions (which is equivalent to (I))

\[ \exists x(Ax\& \sim Bx) \vee (\sim Ax \& Bx) \rightarrow \sim \iota AB\]

Aristotle's own example involves $A$ as 'man' and $B$ as 'living being participating in knowledge'.  The individual $x$ is 'God'. Aristotle states that since $\sim Ax$ and $Bx$ it cannot be the case that $\iota AB$.

The following passage in (133b) is  quite difficult to interpret: (W.A. Pickard translation):\\

\emph{Next, for destructive purposes, see if the property of things that are the same in kind as the subject fails to be always the same in kind as the alleged property: for then neither will what is stated to be the property of the subject in question. Thus (e.g.) inasmuch as a man and a horse are the same in kind, and it is not always a property of a horse to stand by its own initiative, it could not be a property of a man to move by his own initiative; for to stand and to move by his own initiative are the same in kind, because they belong to each of them in so far as each is an 'animal' .}\\

The idea here seems to be that properties of species can vary in a way appropriate to their genus and it in way that is invariant under change of species.  Since standing or moving voluntarily are both appropriate to animal, both properties would have to be shared by animal and horse.

Aristotle's theory of opposition is difficult to fathom from a modern perspective.  Opposition is not classical or intuitionistic negation.  Rather it corresponds, as we saw in some detail in the previous Section, to a second-order function defined on monadic predicate terms. But in fact there are many distinct species of opposition.  $\sim Ax$ does not necessarily imply that $A^\bullet x$.  The following is the fundamental topic for property and contraries:

\[ \iota AB \leftrightarrow \iota A^\bullet B^\bullet \]

In Section 3 we denoted Aristotle's 'modern' notion of term-negation by $T^c$.  The associated topics are not surprising (136a). For instance $\iota AB^c \rightarrow \sim\iota AB$ and $\iota AB \rightarrow \sim \iota AB^c$. 

Chapter 8 of book V deals with degrees of terms as well as with what might be called probabilistic logic. The aspect of degree of a term seems very much to act like an operator on terms which in ESOL would correspond to a second-order monadic function. Thus we can denote 'more','most', 'less' and 'least' $T$ by $T^+, T^\triangle, T^-, T^\bigtriangledown$.  If a term does not admit degree then we assume that $T^+ = T$, etc.  The relevant topics can be condensed into

\[ \iota AB \leftrightarrow \iota A^+ B^+ \leftrightarrow \iota A^- B^- \leftrightarrow \iota A^\triangle B^\triangle \leftrightarrow \iota A^\bigtriangledown B^\bigtriangledown\]

It does not seem that degree has a directly comparative aspect. Rather the degree is taken as relative to some implicit absolute standard. But for probabilistic degree the situation is otherwise.  Let us write $\mathbb{P}\iota AB > \mathbb{P}\iota A'B'$ for 'it is more likely for $B$ to be the property of $A$ then $B'$ to be the property of $A'$ - which is in reality a quaternary second-order predicate in the variables $A,B,A',B'$.  Then Aristotle's rather strange supposition can be expressed directly as

\[ \mathbb{P}\iota AB > \mathbb{P}\iota A'B' \,\&\, \sim \iota A B \rightarrow \sim \iota A'B' \]

In 138a we see variants on the topics corresponding to the cases in which $A = A'$ and $B= B'$.  There is also an analogous (quaternary) relation of equal likelihood: 

\[ \mathbb{P}\iota AB = \mathbb{P}\iota A'B' \,\&\, \sim \iota A B \leftrightarrow \sim \iota A'B'  \]

as well as particular cases in which $A= A'$, etc.

These topics are strange because Aristotle does not state clearly the intrinsic relationship between $A$ and $A'$ or $B$ and $B'$. The following difficult passage seems relevant:\\

\emph{The rule based on things that are in a like relation differs from the rule based on attributes that belong in a like manner, because the former point is secured by analogy, not from reflection on the belonging of any attribute, while the latter is judged by a comparison based on the fact that an attribute belongs.}\\


Book V ends with the curious rejection of property given by superlatives - surely not the same superlative as $T^\triangle$ but rather something like a modifier of a noun, something of the form 'the Xest Y'. Aristotle says basically that such a definite description will not denote uniquely for different states of affairs.

Some genera 'admit degree' and some do not. Take an adjective admitting degree like 'desirable'. Of course this must be specified to being desirable to a certain kind of being.  An adjective admitting degree particularized to a given kind $G$ can be formalized as a pair $(\mathcal{X},\sqsubset)$ where $\mathcal{X}$ represents the set of immediate species of $G$ and $\sqsubset$ expresses comparison regarding the degree of the adjective.  The superlative (if it exists) is unique (152a). For this to be the case we might need the anti-symmetry property $A \sqsubset B\,\&\, B \sqsubset A \rightarrow A = B$. We define a superlative of $G$ relative to $\sqsubset$  to be an element $S \in \mathcal{X}$ such that $\forall A(\gamma^\star AG \rightarrow (A\neq S \rightarrow A \sqsubset S))$.  We write $S = sup^1_\sqsubset G$. There is a problem with the construction 'the + $\langle$ superlative $\rangle$ ', for instance 'the most desirable'. Does this definite article already imply that there is in fact only one superlative (if at all), so that the mentioned argument of Xenocrates is just a case of the topic stated further ahead, that involving two things equal to a third being equal amongst themselves? The only way to make sense of this argument is to cast it in the form:  $S$ is a superlative and $S'$ is a superlative (for $(\mathcal{X},\sqsubset)$) therefore $S = S'$. 

Now Aristotle's furnishes us with an interesting discussion of superlatives applied to pluralities. For instance the 'Spartans' and 'Peloponnesians' are both said to be the bravest of the Greeks.  What do such plural superlative judgments mean ?  How are we to make sense of Aristotle's argument that the Spartans being the bravest and the Peloponnesians being the bravest only means that (the extension of) one term is contained in the other ?

The superlative in question seems to be distinct and  more extensional. We view Spartans and Peloponnesians not as species of 'Greek' but as arbitrary terms having extensions contained in that of 'Greek'. We need a different concept of superlative.  We assume that there is a partial order $<$ defined on the extension of $G$, that is, on individuals $x$ such that $Gx$.  We write as before $A \subset B$ for $\forall x(Ax \rightarrow Bx)$.
Let us posit that $S = sup^2_< G$ for $S \subset G$ with non-empty extension iff

\[ \forall x,y(Gx \&\sim Sx\& Sy \rightarrow x < y) \]

We can now justify (in ESOL) Aristotle's argument about the Spartans and Peloponnesians.

\begin{theorem} If $S = sup^2_< G$ and $T = sup^2_< G$ then $T \subset S$ or $S \subset T$.
\end{theorem}

\begin{proof} Suppose that it is not the case that  $T \subset S$ or $ S \subset T$. Then (since we are assuming that terms $T$ and $S$ have non-empty extension) there is an $a$ such that $Ta$ and such that $\sim Sa$ and a $b$ such that $Sb$ and $\sim Tb$. But then by definition we must have both $a > b$ and $a < b$ and so $a = b$, a contradiction.
	
\end{proof}

We end this Section by illustrating how a passage of chapter 13 of book VI (about 'this and that') might be interpreted (it relates to mereology). Let us define $A+ B = \lambda x(Ax\& Bx)$. Suppose opposition does not satisfy the analogue of De Morgan's laws but rather $(A+B)^\bullet = A^\bullet + B^\bullet$. Indeed it makes sense to say that the opposite of temperate and brave is intemperate and cowardly. Now suppose that we have defined a comprehension operation on individuals which allows us to form the pair $\{x,y\}$ which is also an individual as far as ESOL is concerned. Then we postulate that

\[ A\{x,y\} \leftrightarrow Ax \vee Ay \tag{M1} \]

With these premises we can easily derive Aristotle's paradoxical conclusion that if we define justice to be temperance and courage then considering the pair of men in his example\footnote{one is temperate but cowardly and the other is brave but intemperate.} this pair will be both just and unjust. The elucidation of this paradox would lead to a closer examination of the 'mereology' present in the \emph{Topics}. In chapter 5 of book V there is a discussion of how property relates to the property of the parts of the thing it is a property of. Aristotle considers the special case of things which have like parts such as the sea. In this case the topic seems to be, if we denote the parthood relation by $A\epsilon B$:

\[ \iota AC \& B \epsilon A \rightarrow \iota BC\]

hence 'being the largest body of water' cannot  be the property of the sea. This topic reflects a distinct mereological outlook from (M1).

 \section{The dialectic of the Topics}

  We now discuss in light of the present work the system of argumentation  for which the \emph{Topics} was written. In the Platonic Dialogues (particularly the early\footnote{in a future paper we show the surprising second-order formalizability of the earlier Socratic dialogues as contrasted to the later ones.} and middle ones) as well as in Aristotle there is the fundamental supposition that philosophical knowledge can be gained by finding the right definition of terms (under the proviso that the term in question is acknowledged to exist, that is, have non-empty extension.). This concern is not abandoned in the \emph{Analytics}. The Socratic attitude is that the first step to wisdom consists in testing one's beliefs: even the knowledge of the error of one's beliefs can be considered cathartic and of value. Thus a major part of philosophical activity involved the testing (i.e. attempting a refutation) of proposed definitions of terms - something which Aristotle acknowledges to be easier than finding correct definitions\footnote{A massive amount of arguments were compiled during the eras of the late Academia and first Stoic schools, as testified by accounts of the lives of Chrysippus and Carneades. And in book I of the \emph{Topics} Aristotle seems to imply that  'notebooks'  containing a kind of thesaurus , organized compilation of terms, definitions and argumentative strategies, were employed during debate (or at least to prepare for them). }. In a modern framework how are we to understand (broadly speaking) the dialectical game of the \emph{Topics} ? In ESOL in which we assume we have a finite collection of constants, functions and predicate symbols (which represent the disambiguated vocabulary of a common language.  The non-logical axioms we assume are $T$ representing the topics (or rather a minimal generating set for them ) and $K$ those representing common agreed upon knowledge. The game begins with Alice's proposal of a definition of a certain term $h A D$. Then Bob basically has to try to derive a contradiction $T + K + h AD \vdash \bot$. The way the game was played apparently involved Alice confirming every move made by Bob, as attested in the Socratic dialogues (where it is rare for the 'opponent' not to agree with Socrates' moves). The 'moves' can be $K$ axioms themselves or inferences, either purely logical or based on a topic. A nice example of this dialectical game in practice is given by Plato's \emph{Lysis} although the types of definition considered do not fall within Aristotle's scheme for monadic predicates terms outlined in book VI of the \emph{Topics}.  Further investigation into the interpretative tradition of Boethius (including Cicero's account of the Topics)  and its rich development into the middle ages (cf. the theory of consequentiae and obligationes \cite{pin}) will shed more light on various historical practices of the 'game' of dialectics.

\section{ESOL$_D$}

  There is an alternative version of the existential elimination rules (both first- and second-order) inspired by Hilbert's $\epsilon$-calculus. Just as Peano's $\iota$ operator corresponds to the definite article we can introduce an operator $\tau$ which expresses the indefinite article.  Thus we extend the syntax  so that if we have a formula $\phi(x)$ then $\tau x \phi(x)$ is an individual term and for a formula $\phi(X^{(n)})$ we have that $\tau X \phi(X)$ is a $n$-predicate term. We could also introduce in a similar the way the standard $\iota$ operator.
 We then replace the natural deduction rules $\exists^{i} E$ ($i = 1,2$) with the rules:
 \begin{center}
 	\AxiomC{$ \exists x \Phi (x)$}
 	\RightLabel{$\exists^1 E_{\tau}$}
 	\UnaryInfC{  $\Phi ( \tau x \Phi(x)) $}
 	\DisplayProof 	
 \end{center}
 
 and 
 \begin{center}
 	\AxiomC{$ \exists X \Phi (X)$}
 	\RightLabel{$\exists^2 E_{\tau}$}
 	\UnaryInfC{  $\Phi ( \tau X \Phi(X)) $}
 	\DisplayProof 	
 \end{center}

 The operator $\tau$ and the above rules allow us (once suitably adapted to a bounded form)  to give a formal interpretation of Aristotle's demonstrative technique of \emph{ekthesis} which calls for a generalized indefinite version of the Peano $\iota$ operator. This will be discussed in great detail in Section 8.  If we have that there exists an A, then we can form the term 'an A' so as to be able to postulate in the context of a proof: let a be an A. An example of a natural language expression of the rule above: every man has a father, hence,  every man's father is the father of that man (or: if every man has a father then a father of a man is that man's father).

There is an alternative version of ESOL in the style of dependent type theory in which all quantifiers are bound by some formula. All quantifiers are bound by a formula $\Psi(x)$ and written   $\forall_{\Psi(x)}\Phi(x)$, which is to be interpreted as $\forall x(\Psi(x) \rightarrow \Phi(x))$ and $\exists_{\Psi(x)} \Phi(x)$ which is to be interpreted as $\exists x(\Psi(x) \& \Phi(x))$.
We can reformulate the quantifier rules (with the $\tau$-operator version of the rules for $\exists$) as follows:

 \begin{center}
 	\AxiomC{$\Psi \rightarrow \Phi$}
 	\RightLabel{$\forall_D^1 I$}
 	\UnaryInfC{$\forall_{ \Psi^y_x (x)}\Phi^y_x$}
 	\DisplayProof
 	\AxiomC{$\forall_{\Psi(x)}\Phi$}
\AxiomC{$\Psi^x_t $}
 	\RightLabel{$\forall_D^1 E$}
 	\BinaryInfC{$\Phi^x_t$}
 	\DisplayProof
 \end{center}

 \begin{center}
 	\AxiomC{$\Psi \rightarrow \Phi$}
 	\RightLabel{$\forall_D^2 I$}
 	\UnaryInfC{$\forall_{\Psi^{Y^n}_{X^n}(X^n)} \Phi^{Y^n}_{X^n}$}
 	\DisplayProof
 	\AxiomC{$\forall_{\Psi(X^n)}\Phi$}
\AxiomC{$\Psi^{X^n}_{T^n} $}
 	\RightLabel{$\forall_D^2 E$}
 	\BinaryInfC{$\Phi^{X^n}_{T^n}$}
 	\DisplayProof
 \end{center}

 \begin{center}

 	\AxiomC{$\Psi^x_t \& \Phi^x_t$}
 	\RightLabel{$\exists_D^1 I$}
 	\UnaryInfC{$\exists_{ \Psi(x)}\Phi$}
 	\DisplayProof
	\AxiomC{$ \exists_{\Psi(x)} \Phi (x)$}
 	\RightLabel{$\exists_D^1 E_{\tau}$}
 	\UnaryInfC{  $\Psi( \tau_{\Psi(x)} \Phi(x)) \& \Phi ( \tau_{\Psi(x)} \Phi(x)) $}
 	\DisplayProof 	
 
 \end{center}

 \begin{center}
	\AxiomC{$\Psi^{X^n}_{T_n}\& \Phi^{X^n}_{T^n}$}
 	\RightLabel{$\exists_D^2 I$}
 	\UnaryInfC{$\exists_{ \Psi(X^n)}\Phi$}
 	\DisplayProof
 	\AxiomC{$ \exists_{\Psi(X)} \Phi (X)$}
 	\RightLabel{$\exists_D^2 E_{\tau}$}
 	\UnaryInfC{  $\Psi(\tau_{\Psi(X)}\Phi(X))\&\Phi ( \tau_{\Psi(X)} \Phi(X)) $}
 	\DisplayProof 	
 \end{center}
 
with the suitable provisos for the universal quantifier introduction rules.  In is interesting to express $\forall_D^i I$ for $i = 1,2$ in sequent form:

 \begin{center}
 	\AxiomC{$\Gamma, \Psi \vdash  \Phi$}
 	\RightLabel{$\forall_D^1 I$}
 	\UnaryInfC{$ \Gamma \vdash \forall_{ \Psi^y_x (x)}\Phi^y_x$}
 	\DisplayProof
 
 	\AxiomC{$\Gamma, \Psi \vdash \Phi$}
 	\RightLabel{$\forall_D^2 I$}
 	\UnaryInfC{$\Gamma \vdash \forall_{\Psi^{Y^n}_{X^n}(X^n)} \Phi^{Y^n}_{X^n}$}
 	\DisplayProof
 	
 \end{center}

where $y,Y^n$ are not free in $\Gamma$.

This corresponds exactly to the actual practice of proof and argument. In a certain context $\Gamma$ we assume a hypothesis regarding a generic $x$ in the form: 'assume we have $\Psi(x)$ for some $x$'. Then we show that $\Phi(x)$ and conclude that 'for all $x$ such that $\Psi(x)$ we have $\Phi(x)$'.

We can also consider a variant of $\bot_c$ which is that of mathematical practice. In sequent form it is expressed.

\begin{center}
 
 	\AxiomC{$\Gamma, \Psi, \sim \Phi \vdash \sim \Psi $}
 	\RightLabel{$\bot'_c$}
 	\UnaryInfC{$\Gamma, \Psi \vdash \Phi$}
 
 	\DisplayProof
 \end{center}

With these variations in the syntax and rules of ESOL we obtain the system denoted by ESOL$_D$.  This variant is arguably very close to the mechanisms for logical expression and argumentation in natural language and it will be used in the Sections ahead. In what follows for convenience we write $\forall_{x,y:\psi}\phi(x,y)$ for $\forall_{x:\psi(x)}\forall_{y:\psi(y)} \phi(x,y)$ and similarly for the other types of variables and quantifiers.
 
 \section{Quantifier rules of ESOL$_D$ in Aristotle and Galen}

In this Section for simplicity we will not distinguish between quantification over first and second-order variables. Thus for instance we refer to $\forall_D I$ to refer to both $\forall_D^1 I$ and $\forall_D^2 I$. 
\subsection{$\forall_D E$}

  The theory of the syllogism in the Prior Analytics has been the subject of intense investigation in recent decades, specially the modal syllogism.  But there is an important overlooked point.  A syllogism is  supposed to be a rule of deduction and the collection of the syllogisms considered valid is supposed to be a system of rules which can be deployed to form complex proofs or deductions.  But what Aristotle does in the Prior Analytics is engage in a series of proofs aimed at showing that certain combinations of sentences are valid syllogisms while others are not.  But we can ask: what were the logical rules that Aristotle used in order to prove which syllogisms were valid and which were not ? And to express the syllogisms themselves or things about syllogisms what was the necessary logical or syntactic structure of the sentences employed ?

We can first observe that the most basic traditional syllogism: 'All men are mortal', 'Socrates is a man' therefore 'Socrates is mortal', which is considered an instance of Barbara, is clearly a form of $\forall_D E$, even if unfortunately there are no treatments of singular (like 'Socrates is a man')  or indefinite premises in the extant works of Aristotle. 

But how is the Barbara syllogism expressed ? Let $Term(A)$ be the second-order predicate expressing that a monadic term $A$ is an Aristotelian term and $A\subset B$ be the second-order relation to be interpreted as $B$ is said of $A$. Then we have

\[\forall_{X,Y,Z: Term} \, X \subset Y \& Y \subset Z \rightarrow  X\subset Z \]

 This is patently a universally quantified sentences which involves multiple generality\footnote{ We say that a sentence has multiple generality if it contains at least one relation (for instance a verb with subject and direct object) in which two distinct arguments are governed by quantifiers. A basic example is ‘every man has a father’ which is equivalent to ‘for any man there is some man such that that man is their father’. Here the pronouns which are the arguments (pronouns with anaphora) of the relation ‘father’ are governed by two distinct quantifiers.
} for the relation ‘$\subset$’.  And furthermore it is clear that this syllogism considered as a sentence is meant to be used in conjunction with successive applications of $\forall_D E$ in order to obtain particular valid instances of the syllogism. Let $M$ represent the term 'Man',  $B$ 'Biped' and $A$ 'Animal'. Then

 \begin{center}
 	\AxiomC{$\forall_{X,Y,Z: Term} \, X \subset Y \& Y \subset Z \rightarrow X \subset Z$}
\AxiomC{$Term(M)$}
 	\RightLabel{$\forall_D^2 E$}
 	\BinaryInfC{$\forall_{Y,Z: Term} \, M \subset Y \& Y \subset Z \rightarrow M \subset Z$ }
\AxiomC{$Term(B)$}
	\RightLabel{$\forall_D^2 E$}
 	\BinaryInfC{$\forall_{Z: Term} \, M \subset B \& B \subset Z \rightarrow M \subset Z$ }
\AxiomC{$Term(A)$}
\RightLabel{$\forall_D^2 E$}
 	\BinaryInfC{$ M \subset B \& B \subset A \rightarrow M \subset A$ }

 	\DisplayProof

 \end{center}

  Consider the discussion of the axioms that need to be applied in a geometric proof in the Prior Analytics 41b13-22. One axiom is\\

‘for all quantities A B C and D if A = B and C = D then C subtracted from A = D subtracted from B’ \\

  This leaves now doubt about Aristotle’s consciousness of universally quantified sentences involving multiple generality (in this case the relation is that of equality) nor about the legitimacy of the $\forall_D E$ rule for obtaining instances of the rule.
  In fact Galen's relational syllogisms [Barn] involve precisely this. For instance Galen [Gal. p.46] mentions:\\

‘for all quantities A B C if A is the double of B and B the double of C then A is the quadruple of C’ \\

and Galen’s relational syllogism is precisely the $\forall_D E$ rule by which we can derive (from the premises of '8, 4 and 2 being quantities): \\

‘if 8 is double of 4 and 4 is double of 2 then 8 is the quadruple of 2’\\

for which we can then further apply Modus Ponens ($\rightarrow E$), a rule which was indubitably part of Aristotle’s logical theory: in Prior Analytics 53b12-15  we read: ei gar tou A ontos anagkê to B einai (...) ei oun alêthes esti to A, anagkê to B alêthes einai: if is it necessary that if A is true then B is true and if A is true then necessarily B is true.
 Furthermore in Aristotle’s Topics, as we saw extensively in Sections 4 and 5 (the formalization in ESOL carries over easily to ESOL$_D$)  we can interpret a topic as being a certain kind of universally quantified sentence \cite{prim, slo} (almost always involving multiple generality) which was then instantiated via $\forall_D E$ according to the dialectical need at hand, an interpretation which is in accordance to the tradition expounded in Boethius' \emph{De topicis differentiis}: a topic is precisely a ‘maximal proposition’.  
  From Sections 3-5 we note that the multiple generality present in the (possibly multiply) universally quantified topic sentences can be  logically quite complex and the quantifications inside the sentence can be for instance over the individuals which belong to the extensions of the universally quantified terms.

\subsection{$\forall_D I$}

 Consider the remarkable passage from chapter 3 of Book II of the Topics (tr. W.A. Pickard):\\

\emph{Whereas in establishing a statement we ought to secure a preliminary admission that if it belongs in any case whatever, it belongs universally, supposing this claim to be a plausible one. For it is not enough to discuss a single instance in order to show that an attribute belongs universally; e.g. to argue that if the soul of man be immortal, then every soul is immortal, so that a previous admission must be secured that if any soul whatever be immortal, then every soul is immortal. This is not to be done in every case, but only whenever we are not easily able to quote any single argument applying to all cases in common (euporômen koinon epi pantôn hena logon eipein) , as (e.g.) the geometrician can argue that the triangle has its angles equal to two right angles.}\\ 

  What Aristotle is first saying is that in general we do not have 'induction' from the species to the genus. If (essential) accident C applies to species A of  B it does not follow that the same accident C applies to B.   Thus if the human soul is immortal it does not follow that the soul is immortal (every soul is immortal).   Aristotle says that we have to assume that being an accident of the species implies being an accident of the whole genus.  Aristotle contrasts this situation to another situation in which euporômen koinon epi pantôn hena logon eipein - in which we are able to deduce (in a uniform way) a universal conclusion through a single argument, as the geometrician reasons about property of a triangle.  Thus the geometrician may start with 'Let X be a triangle' and arrive at the conclusion that X has its angles equal to two right angles and thereby view this as a proof that all triangles have their angles equal to two right angles. Thus we are lead to interpret Aristotle's \emph{koinon epi pantôn hena logon eipein} as precisely a statement of rule $\forall_D I$. 

 \begin{center}
 	\AxiomC{$\Gamma, Triangle(x) \vdash  \Phi(x)$}
 	\RightLabel{$\forall_D I$}
 	\UnaryInfC{$ \Gamma \vdash \forall_{ Triangle(x)}\Phi(x)$}
 	\DisplayProof

 \end{center}

  It is also interesting to compare this to Proclus' commentary of the first book of Euclid [proc, p.162] which may be also interpreted as expressing $\forall_D I$:\\

\emph{Furthermore, mathematicians are accustomed to draw what is in a way a double conclusion. For when they have shown something to be true of the given figure, they infer that it is true in general, going from the particular to the universal conclusion. Because they do not make use of the particular qualities of the subjects but draw the angle or the straight line in order to place what is given before our eyes, they consider that what they infer about the given angle or straight line can be identically asserted for every similar case. They pass therefore to the universal conclusion in order that we may not suppose that the result is confined to the particular instance. This procedure is justified, since for the demonstration they use the objects set out in the diagram not as these particular figures, but as figures resembling others of the same sort (...)Suppose the given angle is a right angle. If I used its rightness for my demonstration, I should not be able to infer anything about the whole class of rectilinear angles; but if I make no use of its rightness and consider only its rectilinear character, the proposition will apply equally to all angles with rectilinear sides.}\\

  There are a number of proofs in Aristotle's Physics and Organon which make use of variables in the same way as the mathematician’s ‘let X be a triangle’ and which cannot be interpreted in any other way than as employing the $\forall_D I$ rule (often in combination with $\rightarrow E$) and these rules are furthermore applied to sentences containing multiple generality.  A typical example is Physics 241b34-242a49, the proof that \emph{hapan to kinoumenon upo tinos anagkê kinesthai} – everything that is moved is moved by something. The beginning and end of the proof read (tr. Hardie and Gaye):\\

\emph{ Everything that is in motion must be moved by something. For if it has not the source of its motion in itself it is evident that it is moved by something other than itself, for there must be something else that moves it. If on the other hand it has the source of its motion in itself, let AB be taken to represent that which is in motion essentially of itself and not in virtue of the fact that something belonging to it is in motion (...) a thing must be moved by something if the fact of something else having ceased from its motion causes it to be at rest. Thus, if this is accepted, everything that is in motion must be moved by something. For AB, which has been taken to represent that which is in motion, must be divisible since everything that is in motion is divisible. Let it be divided, then, at G. Now if GB is not in motion, then AB will not be in motion: for if it is, it is clear that AG would be in motion while BG is at rest, and thus AB cannot be in motion essentially and primarily. But ex hypothesi AB is in motion essentially and primarily. Therefore if GB is not in motion AB will be at rest. But we have agreed that that which is at rest if something else is not in motion must be moved by something. Consequently, everything that is in motion must be moved by something.}\\

  The proof is reduced to the special special case: everything that is changed in its own right is changed by something. Here ‘everything’ is to be understood as ‘every physical body’. The proof involves starting with a variable AB: Let AB be taken to represent, etc. Thus Aristotle begins by employing a variable AB and assuming the hypothesis for AB:  ‘P holds of AB’, where P represents ‘changed in its own right’. Then after a detour through a reductio ad absurdum argument Aristotle arrives at the conclusion that ‘some body moves AB’.   The conclusion ‘every body such that P holds of it is moved by some body’ must be seen as following precisely from an application ofthe  $\forall_D I$ rule for the variable AB.

  There are passages in the Physics which can be interpreted as involving multiple applications of $\forall_D I$  rule to derive a conclusion, for instance the theorem in 243a32: the agent of change and that which is changed must be in contact (although there are no variables employed explicitly in the text).

\subsection{$\exists_D I$}

  How does Aristotle show in the Prior Analytics that a certain candidate for a syllogism is not valid ? A common strategy involves presenting instances of three terms for which the premises hold but the conclusion patently does not.  

  It is plausible to assume that Aristotle accepted the rule for converting the negation of universal quantification to existential quantification and the negation of existential quantification to universal quantification for any sentence, not only sentences of a restricted type with only one quantification and without relations. In particular he would accept repeated applications to multiple quantifications.

  Thus if he wanted to show that a syllogism was invalid such as in 26a:\\

‘It is not the case that for all terms A B and C if B is said of all A and C is not said of any B then C is said of all  A’\\

\[ \sim \forall_{X,Y,Z:Term}\, A\subset B \& B \subset C^c \rightarrow A\Subset C) \]

he would use this conversion rule three times to obtain the equivalent sentence:\\

‘There are some terms A B and C such that it is not the case that if B is said of all A and C is not said of any B then C is said of all A’\\

\[  \exists_{X,Y,Z:Term}\,\sim ( A\subset B \& B \subset C^c \rightarrow A\Subset C) \]

  Note that this sentence again has multiple generality.  But an implication is false if the premises are true and the conclusion is false, thus the above sentence is equivalent to:\\

‘There are some terms A B and C such that B is said of all A, C is not said of any B and C is not said of all A’\\

\[  \exists_{X,Y,Z:Term}\,  A\subset B \& B \subset C^c \& \sim( A\Subset C) \]

Aristotle gives us the example:

‘Man is term’,  ‘Animal is a term’, ‘Stone is a term’, ‘Animal is said of all men and stone is not said of any animal and stone is not said of all men’. \\

\[ Term(M)\& Term(A) \& Term(S) \,   M \subset A \& A \subset S^c \& M \subset S^c (or \sim (M \subset S))   \]

 Then with this example it is implied that we can apply the $\exists_D I$ rule three times (for ‘Man’, ‘Animal’ and ‘Stone’) to obtain\\

‘There are some terms A B and C such that B is said of all A, C is not said of any B and C is not said of all A’\\

\[  \exists_{X,Y,Z:Term}\,  A\subset B \& B \subset C^c \& \sim( A\Subset C) \]
 
Hence we can assume that Aristotle not only deploys the $\exists_D I$ rule but does so in the context of multiple generality.

 We will present more evidence for the logically conscious use of the $\exists_D I$ rule in what follows where we discuss the evidence for $\exists_D E$, as  Aristotle frequently  used  these two rules together. We note that it would be worth investigating in a future paper the relationship between $\exists_D I$ and a hypothetical rule which states that from a singular or indefinite sentence (for instance in Boethius’ sense which we consider further ahead) we can deduce the corresponding existentially quantified one: from ‘It is a man’ and ‘It is mortal’ we deduce ‘Some man is mortal’.

If we take ‘It is a man’ then the deduction ‘Something is a man’  is taken to employ implicitly the assumption ‘It is a thing’: we interpret ‘Something’ as quantifying over a large but nevertheless defined universe of things. 

Of interest to the passage between singular or definite sentences to existentially quantified sentences and to the $\exists_D I$ rule in general is that fact that many of the proofs in Euclid’s Elements(cf. [heath]) can be viewed as being based on what is called a constructivist interpretation of existential quantification (see for instance\cite{knorr}). The idea is that to prove for instance that ‘For all X there is some Y is such that P holds between that X and that Y’  amounts to nothing more than producing a construction C which is such that for a given generic A which is X, yields a Y which is in the relation P with X. This can be symbolized as ‘Y = C(X)’ and ‘X and C(X)’ are in the P-relation.  

Thus in this case the $\exists_D I$ rule would be read as the following variant:\\

Given a construction C such that for all A which is X we have that C(X) is a Y which is in a relation P with A then we can deduce that for all X there is some Y which is in the relation P with that X.\\

In ESOL$_D$ this could be written as the sequent rule:

 \begin{center}
 	\AxiomC{$\Gamma, \Psi(x) \vdash \Phi(f(x)),  P(x, f(x))$}
 	\UnaryInfC{$ \Gamma \vdash \forall_{ x: \Psi(x)} \exists_{y: \Phi(y)}P(x,y)$}
 	\DisplayProof

 \end{center}

which clearly relates to Skolemization.

\subsection{$\exists_D E$}

 Consider again the proof of the conversion of universal negation in the Prior Analytics 25a14-18 (tr. Hardie and Gaye): \\

First then take a universal negative with the terms A and B. If no B is A, neither can any A be B. For if some A (say C) were B, it would not be true that no B is A; for C is a B. But if every B is A then some A is B. For if no A were B, then no B could be A. But we assumed that every B is A.\\

 Aristotle's argument (expressed in a rather terse and laconic style) can be read  as follows: \\

Thesis: If A is not said of any B then B is not said of any A\\

1) Assume A is not said of any B.

2) Assume that B is said of some A (Aristotle takes this immediately to be equivalent to it not being the case that B is not said of any A).

3) Let this some A be c. Then c is a A which is also a B.

4) Hence some B is A.

5) But this contradicts the initial hypothesis that A is not said of any B.

6) Hence by Reductio ad Absurdum we conclude that it is not the case that B is said of some A, that is to say, B is not said of any A.

7) Hence the thesis follows from discharging 1).\\

  It is important to  note that  (as generally accepted) Aristotle deploys Reductio ad Absurdim in a logically conscious way: he called this rule hê eis to adunaton apodeixis (see for instance Prior Analytics, 41a21-34).

  Let us look in more detail at what is going on at 3). We start with ‘some A is B’.  Aristotle is saying that c is ‘a A which is B’.  But then  ‘a A which is B’ is B, that is, ‘c is B’, and also ‘c is  A’. But then he deduces immediately that ‘Some B is A’. Thus at 3) it is plausible to interpret the argument in 3) and 4) as involving first the $\exists_D E$ rule (where ‘a A which is B’ is conveniently renamed c):\\

From ‘some A is B’ we conclude ‘A A which is B is B’ (i.e. c is B) and ‘An A which is B’ is an A’ (c is A)\\

followed by the application of $\exists_D I$:\\

From ‘c is B’ and ‘c is A’ (note we switched the order) we conclude ‘Some B is A’.\\

 Or writing out 2-4 of this proof in ESOL$_D$.\\

2. $\exists_{x: A(x)} B(x)$

3.  $A(c) \& B(c)$  ($\exists_D E$, 2, for  $c = \tau_{x:A(x)}B(x)$ )

4. $\exists_{x: B(x)}A(x)$ ($\exists_D I$, 3)\\

  In 28a24-27 of the Prior Analytics Aristotle gives two proofs of a certain syllogism in the third figure, the second proof employs a process called \emph{tô(i) ekthesthai poiein}, by \emph{ekthesis}.  We then confirm that this process corresponds to the deployment of $\exists_D E$:\\\

  If they are universal, whenever both P and R belong to S, it follows that P will necessarily belong to some R. For, since the affirmative statement is convertible, S will belong to some R: consequently since P belongs to all S, and S to some R, P must belong to some R: for a syllogism in the first figure is produced. It is possible to demonstrate this also per impossible and by exposition (ekthesis). For if both P and R belong to all S, should one of the Ss, e.g. N, be taken, both P and R will belong to this, and thus P will belong to some R. \\

  Recall that for Aristotle universal predication has existential import so we can assume that ‘some thing is a S’. Aristotle's second proof of this third figure figure (by exposition or ekthesis) can be read in detail as follows: \\

  Assume that ‘all S is P’ and that ‘all S is Q’. We have then that there exists an S (universal predication has existential import, which we read as ‘Some thing is S’).  Let n be a thing that is S. Then by the hypotheses and using $\forall_D E$ we conclude that ‘n is Q’ and likewise also that ‘n is P’ . Thus by $\exists_D I$  ‘some P is Q’.\\

In ESOL$_D$:\\

1. $\forall_{x:S(x)} P(x)$ \& $\forall_{x:S(x)} Q(x)$ (Hyp)

2. $\exists_{x: thing(x)} S(x)$ (Aristotelian axiom)

3. $S(c)$ ($\exists_D E$,2, for $c = \tau_{x: thing(x)} S(x)$, after using $\&E_R$)

4. $P(c)$ ($\forall_D E$ 1,3, after applying $\& E_L$)

5. $Q(c)$ ($\forall_D E$ 1,3, after applying $\& E_R$)

6. $P(c) \& Q(c)$ ($\& I$)

7. $\exists_{x: P(x)} Q(x)$  ($\exists_D I$, 6)

8 $\forall_{x:S(x)} P(x)$ \& $\forall_{x:S(x)} Q(x) \rightarrow \exists_{x: P(x)} Q(x)$  ($\rightarrow I$, 1, 7)\\

  This proof besides presenting strong evidence that proof by \emph{ekthesis} involved precisely using the $\exists_D E$ rule (for Alexander of Aphrodisias ekthesis involved an individual which was 'set out' from the extension of some term ) it also presents additional evidence for Aristotle’s use of $\exists_D I$  in the step 6 -7: that the passage from indefinite (or singular) to existentially quantified sentences was an integral part of his logical theory.\\

  Regarding the ‘singular’ sentences occurring in the proof  it is regrettable that (unlike in Boethius) propositions such as 'Socrates is mortal' are not explicitly discussed in the extant works of Aristotle, specifically with regards to the version of the third-figure: 'Socrates is mortal' and 'Socrates is a man' hence 'Some man is mortal'.

  The above proof does not involve multiple generality, but examples which we can interpret as involving the deployment of $\exists_D E$ to sentences of multiple generality are indeed present in the works of Aristotle.

  We will focus on sentences involving multiple generality of the form\\
 
‘for all X some Y is such that the X and the Y are in relation R’\\

\[\forall_{x: X(x)} \exists_{y: Y(y)} R(x,y) \tag{1}\]

for example ‘for every man there is some man that is that man’s father’.\\

  If we change ‘there exists’ to ‘there exists a unique’ then we get the modern concept of ‘function’.
The use of the $\exists_D E$  rule for sentences of the above form is found in the Physics. 

  Typically we have a variable Z (to be used for $\forall_D I$ and introduced as ‘Let Z be such that...’) and we first use $\forall_D E$ to instantiate the sentence to\\

‘some Y is such that Z and Y are in relation R’\\

and then the application of the $\exists_D E$  rule involves the term ‘a Y such that Z and Y are in relation R’  and yields\\

‘Z and a Y such that Z and Y are in relation R are in relation R’   and ‘a Y such that Z and Y are in relation R is a Y’.\\

  In this case our term ‘a Y such that Z and Y are in relation R’   depends explicitly on Z and it would not be very clear to represent it simply as a constant N; rather is would be better to indicate its dependency by using  a functional notation like N(Z).  The concept of function corresponds to a form of genitive: N(Z) is Z’s Y such that Y is in relation R to it.\\

In terms of ESOL$_D$ we assume the hypothesis $X(z)$ and use $\forall_D E$ and (1) to obtain
\[\exists_{y: Y(y)} R(z,y)\]

Then we can form the term $\tau_{y: Y(y)} R(z,y)$ which has free variable $z$ and thus is a function $n(z) = \tau_{y: Y(y)} R(z,y)$ for which $R(z,n(z))$.

  These logical forms abound in the Physics. For example we have ‘every motion has a given time taken’ and therefrom we find ‘let A be the time taken by motion M’.  Or more generally ‘for every body and every motion of that body there is a time taken by the motion of that body’ and using this we find expressions of the form ‘let T be the time taken by motion M of body Z’. In the associated proofs in the Physics we thus find the use of $\exists_D E$.

  One of the questions which present themselves regarding the Topics is how the definition of relations were carried out.  While there is evidence that Aristotle did possess such a theory (cf. 145A14-15, 142a26), the details are lacking. To end this Section we engage in a speculative reconstruction of how Aristotle would have defined a very basic relation (that of a man being the grandfather of a certain person) and given such a definition how basic reasoning would have been carried out by Aristotle using $\exists_D E$ (or proof by \emph{ekthesis}).

  There is little doubt that Aristotle would have accepted the definition:\\

A being the paternal grandfather of B is A being the father of C and C the father of B for some C\\

as well as the axiom\\

Every man has a father. (H)\\

  We will now attempt to reconstruct how Aristotle would have proven that\\

Every man has a grandfather. \\

  A reconstruction of his argument would be something like:\\

Let A be a man. Then by (H) A has a father. Let X be A's father. Now again by (H) X has a father. Let B be X's father. Then B is the father of X ($\exists_D E$) and X is the father of A ($\exists_D E$). Thus there is some C such that B is the father of C and C is the father of A ($\exists_D I$). So by definition, B is the paternal grandfather of A. So A has a paternal grandfather. Thus every man has a paternal grandfather ($\forall_D I$).\\

Thus we see how $\exists_D E$ could have been deployed in conjunction with other rules in cases of multiple generality.  For simplicity we did not include the additional assumption that being the father of man implies being a man.

 \section{Multiple Generality in Boethius}
In this Section we  look at certain passages from Boethius which can be construed as presenting a formal theory of the logical syntax of sentences which we further argue to include multiple generality. We assume that the passages in question represent a transmission of much more ancient material (including that of the Peripatetic schools).
  Let us consider passages in Book I 1173D-1176D of Boethius' De topicis differentiis, which are of particular interest. We refer directly to the Latin text \cite{boe}.
  In 1174D we have  a classification of propositions which includes the indefinite and singular type of proposition (the detailed treatment of which is so conspicuously absent in the extant works of Aristotle). The indefinite type of proposition is given by Boethius' example 'Homo iustus est' – ‘a man is just’.  We will see further ahead the crucial role this type of proposition plays in quantified sentences and  multiple generality.

  Boethius would accept noun phrases and in particular definite descriptions as subjects. There is evidence for this 1175D:...\emph{partes quas terminos dicimus, non solum in nominibus uerum in orationibus inueniantur} - parts of which are called terms which are found to be not only nouns but also phrases.

  An interesting aspect of Book I is Boethius' account of the logical and syntactic structure of propositions. It is clear that Boethius admits what are called 'simple' propositions, involving the standard subject-predicate term relation - and these simple propositions can be universal, particular, indefinite and singular. The terms occurring can themselves be phrases (orationes) - and we already find this in an example of the 'noun' - 'verb' construction given in Plato's Sophist: \emph{Theaetetus, with whom I am now speaking, is flying}.  But to Boethius the simple proposition is only one type of proposition. The other kind being the conditional proposition which is not simple but complex, being built up from two propositions.

 There are two major questions which naturally arise. The first is: is the classification into universal, particular, indefinite and singular valid for propositions in general (thus also for conditional propositions) or only for simple propositions ? The second one is: is the definition of conditional proposition inductive, that is, can the 'parts' of a conditional proposition be any proposition (in particular be themselves a conditional proposition) or are they forced to be simple propositions ? 

 There is no doubt that the structure of the text – a classification by division - allows us to give an affirmative answer to the first question. After the division into universal, particular, indefinite and singular, Boethius proceeds to divide each of these kinds: 1175A  Harum uero alias praedicatiuas alias conditionales uocamus. It is of utmost interest to examine the logical and grammatical structure of the examples adduced. But here we confine ourselves to the universal conditional proposition. Note that  'all men are mortal' is expressed ‘for all beings if that being is a man then that being is mortal’ or more succinctly ‘if it is a man it is mortal’. Boethius would say that the antecedent in the simple proposition 'it is a man' and the consequent the simple proposition 'it is mortal'. What kind are such simple propositions? They can only be indefinite. We thus see the key importance of indefinite propositions for the expression of their subjects via personal pronouns with anaphora, the way they link up subjects from different propositions of a conditional. 

  It is Boethius' example of a conditional involving indefinite (or  perhaps singular ) simple propositions that allows us to propose an interpretation of sentences expressing multiple generality such as: 'if something happens then something follows it according to a rule'. Boethius would analyze this as having the antecedent ' it happens' and consequent 'something follows it according to a rule'.  As for the second question about the definition of conditional being inductive, the answer is indubitably affirmative 1176B: \emph{Harum quoque aliae sunt simplices conditionales aliae coniunctae. Simplices sunt quae praedicatiuas habent propositiones in partibus (...) Coniunctarum uero multiplex differentia, de quibus in his uoluminibus diligentissime perstrinximus, qua de hypotheticis composuimus syllogismis.}

  The sentence in 1176A 'Conditionalium uero propositionum, quas Graeci hypotheticas uocant partes, sunt simplices propositiones' is difficult to translate and interpret consistently. Is Boethius saying that simple conditionals are called hypothetical conditionals by the Greeks ? Or that the Greeks call the parts of conditionals hypothetical parts ? But then 'sunt simplices propositiones' does not agree with the discussion on conjunct conditionals which are seemingly also called hypothetical syllogisms.

  Only the species of simple conditionals are required to have simple propositions as their parts (their antecedent and consequent) but the same is not true for conjunct conditionals. 

  Boethius in the passage cited above gives us testimony of the existence of various detailed works on conjunct conditionals : \emph{Coniunctarum uero multiplex differentia, de quibus in his uoluminibus diligentissime perstrinximus, qua de hypotheticis composuimus syllogismis}.

  Where does multiple generality come in ? Recall that the parts of a universal conditional proposition can be any kind of proposition including those with phrase terms. And that it seems that the propositions forming part of the  conditional are considered in a special grammatical way via pronouns and anaphora. Consider the sentence: 'if something happens then something follows it according to a rule'. Boethius would analyze this as having the antecedent ' it happens' and consequent 'something follows it according to a rule'. If Boethius accepts in 1176A the phrase predicate ‘investigates the essence of philosophy’ (instantiation of a relation term by a singular term) why would he not, since he admits phrases as predicate or subject terms, accept the predicate ‘being followed by something according to a rule’ ? Further light on this matter might be obtained by analyzing the possible cases of the parts of universal conditional propositions and how the entanglement of pronouns (by anaphora) worked in such cases. 

 We now present  a reconstruction of a fragment of formal logical syntax of the account of propositions given by Boethius in Book I of \emph{De Topicis Differentiis}.  This fragment only involves simple conditionals (those involving an antecedent and consequent consisting of simple propositions) and does not include ‘conjunct’ conditionals concerning which, as we saw, Boethius claimed to have written about in detail in several volumes. For simplicity we do not treat negation (or assume that we are considering term-level negation only).  The most interesting aspect of this logic is that it includes indefinite and singular propositions,  that conditionals can themselves (irrespective of their components) be either universal, particular, indefinite or singular (and a major challenge is to give account of the types of non-universal conditionals) and that entire phrases can function as terms (as subject or predicates). 

 Let us start with the simple proposition which has the form\\

S is P\\

where S and P are ‘terms’.   We consider that terms have a ‘valency’  which determines how many arguments they can be applied to (including the subject).  Thus in the proposition ‘S is P’  the term P must have valency 1.  But the term P (as well as S) can be a ‘phrase’ (oratio).  The simplest type of such term (present in the example given by Boethius himself) is given by a verb-term V followed by a direct object term O. The verb-term V must evidently have valency 2 and thus the complex term VO have valency 1 in order to be applied to the subject term S.

 But what about the types universal, particular, indefinite and singular of simple propositions. The first two evidently correspond to ‘all S is P’, ‘some S is P’ and the singular proposition ‘S is P’ where S denotes a single object such as Socrates or the sentence ‘1 + 1 = 2’.  We can also write singular propositions in the modern form P(S). 

 Indefinite propositions seem to play a central role in conditional propositions and we argue that they correspond to certain extent to our modern concept of expression with a free variable, the free variable in this case being a third person pronoun (which we assume to be ‘it’).  Thus the indefinite simple proposition is ‘it is P’ which we can think of as P(x). 

 The first hypothesis in our reconstruction involves assuming (and this seems plausible) that Boethius allowed for abstraction on propositions which yielded valency 0,  1 and even 2 terms.

 That is to say, from a proposition ‘some S is P’   he would allow us to derive the following phrase terms (the cases for ‘all’ and singular propositions are similar):\\

(being) the case that that ‘some S is P’ (valency 0)

(being) such that some of it is P (valency 1) 

(being) such that it holds for some S (valency 1)\\

For the case ‘some S is VO’ we could further derive the phrase term:\\
 
(being) such that some S is V of it (valency 1)\\

and analogously we could derive from ‘all S is VO’ the phrase term:\\

(being) such that all  S is V of it (valency 1)\\

 Thus Boethius could make use of this last phrase predicate to obtain: \\

all S is VO for some O = for some O all S is V of that O\\

  Consider the example ‘every man knows some things’.  Then the proposition ‘some things are known by Socrates’  and then apply the abstraction above to obtain the 1-valent term ‘(being) such that some things are known by it’. We then obtain\\

all man is such that some things are known by them\\

  Also we note that from a simple proposition ‘S is VO’  it is likely that Boethius would allows us to derive the indefinite proposition  ‘S is V of it’.\\

  Let us now consider simple conditionals.  The form of the universal simple condition is clear. It is of the form\\

if it is A then it is B      or,    (for all thing) if it is A then it is B.\\

which Boethius considers as having components the indefinite propositions ‘it is A’ and ‘it is B’ which in modern form would be A(x) and B(x) or A(it) and B(it).  What is fundamental here is that the two variables are linked grammatically: this corresponds to the anaphora of the  pronoun ‘it’ as discussed by Bobzien and Shogry\cite{mult} for this case of Stoic Logic (but for a different set of pronouns corresponding to  classical Greek tis and ekeinos).  We could also express our multiple generality example above as the universal simple conditional\\

 if it is a man then some things are known by them\\

 Our challenge now is to give an account of particular, indefinite and singular simple conditionals.  We propose for particular conditionals\\

 something is such that if it is A then it is B\\

and for indefinite conditionals (note that these sentences, being themselves indefinite, could be used to further form ‘conjunct’ conditionals)\\

 it is such that if it is A then it is B\\

 With this in mind, we can speculate that an example of a universal conjunct conditional would have been:\\

 (for any thing) if it is such that if it is A then it is B then it is C\\

The most difficult case is that of  singular conditionals and indeed it is difficult to see were the standard modern conditional ‘if P then Q’ comes in. 

Our guess for an example of a singular simple conditional would be:\\

if Socrates is A then Socrates is B\\

but why not also\\

if Socrates is A then Plato is B\\

and indeed the logical form suggests that singular simple conditionals could in general include antecedent and consequents which can either be universal, particular or singular (and thus correspond to standard propositional implication).

From what we have seen further investigation is required related to how universal quantification, universal bounded quantification (as in ESOL$_D$) and universal simple conditionals are related 
not only in Boethius but in Aristotle himself. It is possible that the above interpretation of Boethian conditionals should be replaced by a bounded version. Instead of 'if it is a man then it is mortal' we should be more through in specified the domain of 'it' and express this more fully as 'for all beings if such a being is a man then that being is mortal' and for the particular  the things over which 'something' ranges must be specified: 'there is a C such that if it is  A then it is B'. 

It is interesting to note that in some versions of dependent type theory implication is subsumed as a special case of what amounts to bounded universal quantification.

 \section{Conclusion}

We have discussed in some detail  how ESOL might be used to formalize the \emph{Topics}. If ESOL is the most natural framework to read the \emph{Organon} as a whole then this sheds new light on the development of Western logic, suggesting a much greater formal and conceptual continuity between ancient logic and Frege's system of the \emph{Begriffsschrift}, usually equated to second-order logic with a term-forming extension operator.

ESOL (in its standard form or with bounded quantifiers\footnote{and Dependent Type Theory, also called Constructive Type Theory.}) does not in itself  formalize extensional logic. It  need not even come with equality predicates (there must be more than one: equality must be divided between a binary first-order predicate for individual equality and  second-order binary predicates for first-order $n$-ary predicates). Even allowing a standard Leibnizian axiomatics of equality, Frege still had to introduce extensionalism as something derived from an ad hoc axiom (Law V). It follows from this law that two predicates are equal iff they apply to the same elements (they have the same 'extension' or 'course-of-values').

In the post-Fregean history of philosophical logic there were schools which posited extensionalism as a dogma and strove furthermore to ban intensionality and modality from logic altogether or to explain away these in an extensionalist anner.  Other logicians, mathematicians and philosophers (including Husserl himself) viewed extensionalism and Frege's Law V as the ultimate source of Russell's paradox and a main source of other paradoxes and confusions in the philosophy of logic (this is explored in particular in \cite{hill}).  We agree entirely with the thesis of \cite{lob} that one can view the remarkable developments in 20th century logic (and I would include homotopy type theory) as a progressive liberation from the shackles of extensionalism and (revidindicating Husserl's program) a gaining of awareness of the essential, inescapable, intensional and modal aspects of  logical systems including classical logic itself (further evidenced by the theory of control operators and the $\lambda\mu$ -calculus). I would add that the current view would be that extensionalism does not flow from the nature of logic itself but is something added and posited for the needs of certain specific 'regional ontologies' .  A major thesis of our paper is that in Aristotle we find what amounts to a rich purely intensional extended second-order logic and at the same time a more specialized development of an extensionalist logic (including with a possible worlds semantics) in the Prior Analytics.  Thus both Frege and post-Fregean developments are likewise anticipated.

A major aspect of the present work is the thesis that Aristotle's theory of logic in the \emph{Topics} is much more general,  richer and philosophically insightful than the subsequent more specialized syllogistic theory of the \emph{Prior Analytics} with the possible exception of some technical aspects of the modal syllogistic.  For monadic first-order predicates, our arguments amount to showing that a major difference lies in that second-order relations used in the Prior Analytics are extensional while those used in the Topics are more fine-grained and intensional.  The transition between the second-order relations of the Topics (such as $\gamma$, $\iota$) and that of the Analytics is one of a flattening and collapse wherein the modal differentiations, the different ways of a predicate belonging to another predicate, are discarded in favor of a purely extensional criterion of predication ($\kappa$).  But the extended second-order logic of the Topics is essentially an intensional logic. 
 
Having an intensional extended second-order logic allows us to formalize a potentially unlimited variety of de dicto modalities by using second-order relations $\rho AB$ applied to predicates $A$ and $B$, expressing how $B$ is related to $A$ independently of the extensions of these predicates. De re modalities  can likewise be formalized through second-order functions as briefly discussed in Section 3.2. In this discussion we point out how this alternative approach can also interpret our AML system we developed for the modal syllogistic (which however used a variant of standard extensionalist possible worlds semantics). It may be that the a non-extensional extended second-order approach to de re modalities for monadic predicates could help refine this interpretation.

Returning to de dicto modalities, we explored a little in Section 5 the richness of these modes present in Aristotle's Topics. Aristotle defines a second-order quaternary relation which cannot be interpreted except as a precursor of probabilistic logic. We wrote the relation as $\mathbb{P}_\iota AB > \mathbb{P}_\iota A'B'$ but it is understood that this is a single quaternary relation $\mathbb{P}ABA'B'$ which expresses that $B$ is more likely to be the property of $A$ than $B'$ is to be the property of $A'$.  Aristotle then presents topoi which express several interresting axioms concerning this relation.

Also (as we saw when we discussed the problemsf of the identity of a term and its definition or  the treatment of identity in Book VII of the Topics) a major avenue of exploration  concerns the different  fine-grained (in particular qualitative) concepts of equality and identity (which include the very particular and specialized notion of extensional equality).

In our investigations of the complex textual  and interpretative problems surrounding the relations $\delta, \Delta$ and accident $\sigma$ (specially in Section 4.1) the variety of kinds of intensional modes of predication becomes apparent as well as the connection to the standard modalities (i.e. those that appear in the modal syllogistic and which are already slightly richer than the standard modern ones).  Although we cannot offer yet a definite account, it appears to be a kind of ontological gradation starting from the center and progressing to the periphery of a being or concept. At the heart is the relation of definition expressing core essence, then follows essential and inseparable accidents, then property (which does not express essence and yet is inseparable) and finally different kinds of separable 'accident' which are defined by contingency, potentiality and genericity (something being the case most of the time).  There is also a relation of having a predicate 'naturally' or 'by nature'. It would be very interesting to consider these relations in light of Husserl's mereology proposed in the Third Logical Investigation.

 Regarding Section 6 we can of course question if this dialectical game could lead to actual progress in philosophy or to certain or plausible knowledge about reality. In particular the exact scope of $K$ seems to be vague and uncertain; the assumption of  the possession of a  perfectly disambiguated language is itself implausible (although it should be said that the \emph{Topics} is acutely aware of this problem). However we should not go to the other extreme of discarding this kind of philosophical practice. Rather what was required is a preliminary  epistemic, semantic and logical analysis of our concepts and language - and such concerns are in no way foreign to the mature Aristotle: things better known to us are not necessarily better known in themselves\footnote{as we saw in Section 6, the extra-logic second-order relation of 'being better known' plays an important role in the \emph{Topics}: if $A$ is a property of $B$ then it has to be better known than $B$. In book I,1 Aristotle says that a syllogism is a 'proof' when the premises are 'true and primitive' or epistemically depending on 'true and primitive' premises.  The adjective 'alethê'  is also translated as 'evident'. Aristotle then clarifies that evident and primitive is what is known through itself and not through another'. All 'scientific principles' should be self-evident.  Other propositions are 'endoksa', translated sometimes as 'probable'.}.  The theory of the \emph{Categories} is already assumed in the \emph{Topics}. We can ask if Aristotle envisioned  the project of a systematic and complete classification of all terms (the entire semantic universe) according to the theory of definition discussed in this paper, starting from the supreme genera. Such a project (which contrasts in certain aspects with that of Leibniz' characteristica)  bears some similarity to the goal of Roget's  famous Thesaurus or  of modern lexicological projects.


  As proposed in Section 4.3, in order to gain a greater understanding of the logical, semantic and grammatical content of the \emph{Topics} it would be  useful to compile a list of all the definitions (both accepted and rejected) mentioned in this work and to pay close attention to their grammar and syntax and how this relates to their semantic and logical structure. We also need to clarify in greater detail the concept of predication deployed throughout the work and how it related to the concrete forms of predication given by the predicabilia.

  We Section 8  have argued for and provided evidence that Aristotle was in possession of a theoretical logic completely adequate for reasoning about quantified sentenced containing multiple generality (together with the basic rules for implication and reductio ad absurdum). This thesis was further corroborated with material from Galen and Boethius.  This work can be seen to parallel in the spirit of its historical approach to the work of Bobzien\cite{bob1} and Bobzien and Shogry\cite{mult} on the sophistication of the deductive system of Stoic propositional logic and the adequacy of this logic to express multiple generality.  While much remains to be studied and explored regarding quantifier reasoning and multiple generality in medieval logic the conclusion which seems to suggest itself is that our standard view of the history of western logic needs to be somewhat revised so that a greater continuity is acknowledged between ancient logic and the development of quantifier logic by Frege and Peirce as well as the calculus of relations of Schröder.  A remarkable consequence is that if our interpretation and reconstruction of Boethius’ logic given in this paper be correct then it furnishes some vital clues to better understand Kant’s conception of logic as expounded in the Critique of Pure Reason and frees Kant to a certain extent (in so far as he would have been acquainted with Stoic logic and the writings of Boethius) from logical naivety or the imprisonment in first-order monadic logic,  a logic unsuitable for mathematics and physics as well as philosophy.
  	
  	\section*{Appendix 1:  key topics of book IV}
  	
  	Here we list the ESOL interpretations of several important logical topics in book IV. To this end we make use of the list compiled by Owen in \cite{ow}[pp.691-694]. We divide the list into the corresponding chapters of book IV. We introduce an additional predicate  $\mathcal{C}AB$ to indicate $A$ and $B$ belong to the same category. We take $\gamma AB$ to encompass  (as already mentioned)  what Aristotle calls '$A$ participating of $B$' and '$A$ receiving the definition of $B$' (121a). We use extensional containment $A\subset B$ for unspecified 'predication'.\\

  	Chapter 1.
  	
  	\begin{enumerate}
  	\item  $Ax \& \sim Bx \rightarrow  \sim \gamma AB$ or $\gamma CA \& \sim \gamma CB \rightarrow \sim \gamma AB$ 
  	\item  $\sigma A B \rightarrow \sim \gamma AB$ and  $\sigma' A B \rightarrow \sim \gamma AB$
  	\item $\sim \mathcal{C} AB \rightarrow \sim \gamma AB$
  	\item $\gamma B A \rightarrow \sim \gamma AB$
  	\item $\exists x(Ax\& \sim  Bx) \rightarrow \sim \gamma AB$ (same as 1)
  	\item $\sim \gamma^\star B A  \& \forall C(\gamma^* C A \rightarrow \sim \gamma B C) \rightarrow \sim \gamma B A$
  	\item $ B\subseteq A \rightarrow \sim \gamma AB$
  	\item $Ax \& Ay \& (\sim Bx \vee \sim By) \rightarrow \neg \gamma AB$  	\end{enumerate}
  	\vspace{10mm}
  	Chapter 2\\
  	
  	\begin{enumerate}
  		\item $ \gamma AC \& \sim (\gamma BC \vee \gamma CB) \rightarrow  \sim \gamma AB$ and in 4: $\gamma^iAx \& \gamma^i Bx \rightarrow (\gamma AB \vee \gamma BA)$
  		\item $\gamma BC \& \sim \gamma AC \rightarrow \sim \gamma AB$
  		\item $\gamma BC \& \gamma CA \rightarrow \sim  \gamma AB$ and also Ch.1,4
  		\item $ Ax \& \sim Bx \rightarrow \sim \gamma AB$ and $Ax \& \gamma BC \& \sim Cx \rightarrow \sim \gamma AB$
  		\item $hBD \& A\not\subset D \rightarrow \sim \gamma AB$ and $hBD \& \gamma CB \& C \not\subset D \rightarrow \sim \gamma AB$
  		\item $\delta DA \rightarrow \sim \gamma DA \& \sim \gamma AD$ and $\Delta AD  \rightarrow \sim \gamma DA \& \sim \gamma AD$
  		\item $\gamma BA \rightarrow \sim \gamma AB$ (cf. Ch.1,4)
  		\item included in 6. (122b) $\Delta AB \rightarrow A\subseteq B$.
  		\item included in 6. (123a) $\delta AB \rightarrow A \subseteq B$ (?)
  		\item $\gamma AB \rightarrow \sim \Delta AB$
  		\item $\forall D(\delta BD \rightarrow A \not\subset D) \rightarrow \sim \gamma AB$
  		\item (not logical)
    	\item $\exists x(Ax \&((\lozenge B^c)x \vee (\delta DB \& Dx \& (\lozenge D^c)x)) ) \rightarrow \sim \gamma AB$
  	\end{enumerate}
  	
  	Chapter 3\\
  	
  	\begin{enumerate}
  		\item $B\neq B^\bullet \& \gamma AB^\bullet \rightarrow \sim \gamma AB$
  		\item $\sim \gamma AA$
  		\item $\sim ( \exists C,D(\gamma^\star CB \& \gamma^\star DB \& C\neq D))\rightarrow \sim \gamma AC$
  		\item (not logical)
  		\item  $B = B^\bullet \& \gamma A B \rightarrow \gamma A^\bullet B$ and $B \neq B^\bullet \& \gamma A B \rightarrow A \neq A^\bullet$ and $B = B^\bullet \& \gamma A^\bullet B \rightarrow \gamma A B$ and $B \neq B^\bullet \& \gamma A^\bullet B^\bullet \rightarrow \gamma A B$  
  		\item (we cannot find this in the text)
  		\item  $\gamma A B \rightarrow \gamma B^c A^c$
  	\end{enumerate}
  	
  	Chapter 4\\
  	
  9. 	  $\sim\forall x Sx (\iota y R^{\ast}xy) \rightarrow  \sim \gamma_2 RS$ (124b)\footnote{if $R$ is one-to-one.}\\

								\end{document}